\documentclass[12pt]{article}
\usepackage{amssymb,latexsym,amsmath}
\def\C{\mathbb C}
\def\R{\mathbb R}
\def\N{\mathbb N}
\def\Z{\mathbb Z}

\numberwithin{equation}{section}

\newtheorem{thm}{Theorem}[section]
\newtheorem{lem}{Lemma}[section]
\newtheorem{con}{Conjecture}[section]

\newtheorem{prop}{Proposition}[section]
\newtheorem{cor}{Corollary}[section]
\begin{document}
\sffamily

\title{Non-real zeros of derivatives of  meromorphic functions }
\author{J.K. Langley}
\maketitle

\begin{abstract}
A number of results are proved concerning non-real zeros of derivatives of real and strictly non-real meromorphic functions in the plane.
\end{abstract}

\section{Introduction}

If $f$ is a non-constant meromorphic function in the plane then so is the function
$$
g(z) = \widetilde f (z) = \overline{ f( \bar z)} .
$$
Here $f$ is called real if $g=f$, and strictly non-real if $g/f$ is non-constant.
If $f$ and $g = \widetilde f $  have zeros and poles at the same points with the same multiplicities, 
which will certainly be the case if all zeros and poles of $f$ are real, then $g/f$ has no zeros 
and poles and has modulus $1$ on $\R$, and so $\widetilde f = e^{ih} f $,
where $h$ is a real entire function.

There has been extensive research into the existence of non-real zeros of derivatives of real entire or meromorphic functions
\cite{BEpolya,BEL,CCS2,HSW,HelW3,KiKim,lajda,Laams09,Lawiman13,rossireal,SS}, 
but rather less in the strictly non-real case. Meromorphic functions  which, together with all their derivatives, have
only real zeros were classified in \cite{Hin1,Hin2,Hin3}. The only
other general result treating the strictly non-real case appears to be the following  \cite[Theorem 1]{HSWsnr}. 

\begin{thm}[\cite{HSWsnr}]
 \label{thmsnr}
Let $f$ be a strictly non-real meromorphic function in the plane with only real poles, such that
$f$, $f'$ and $f''$ have only real zeros. Then $f$ has one of the following forms:
\begin{eqnarray*}
(I) \quad f(z) &=& Ae^{Bz}  \quad  ; \\
(II) \quad f(z) &=& A \left(e^{i(cz+d)} - 1 \right)  \quad  ; \\
(III) \quad f(z) &=& A \exp( \exp( i(cz+d) ))  \quad  ; \\
(IV) \quad f(z) &=& A \exp \left[ K( i(cz+d) - \exp( i(cz+d) ) ) \right]  \quad  ; \\
(V) \quad f(z) &=&  \frac{A \exp [ - 2i (cz+d) - 2 \exp( 2i(cz+d)) ] }{\sin^2 (cz+d)}    \quad  ; \\
(VI) \quad f(z) &=&  \frac{A}{e^{i(cz+d)}-1} \quad .
\end{eqnarray*}
Here $A, B \in \C$, while $c, d$ and $K$ are real with $K \leq -1/4$. 
\end{thm}

In the last example (VI) it is easy to verify that $f$ is strictly non-real but $f'$ is not, while 
$f$ and $g = \widetilde f$ have no zeros, and the same poles, 
and $f^{(m)}$ and $g^{(m)}$ have the same zeros for all $m \geq 1$; moreover, 
$f'$ has no zeros, and $f''$ has only real zeros,
but if $m \geq 3$ then $f^{(m)}$ has infinitely many non-real zeros, by \cite[Lemma 3.1]{Laams09}.
The following theorem will be proved, and uses standard terminology from \cite{Hay2}. 

\begin{thm}
\label{thm1}
Let $f$ be a strictly non-real meromorphic function in the plane, and assume that:\\
(i) $f$ has finitely many zeros;\\
(ii) $f$ has finitely many non-real poles;\\
(iii) $f^{(m)}$ has finitely many non-real zeros for some $m \geq 2$.

Then the Nevanlinna characteristic of $f'/f$ satisfies 
\begin{equation}
 \label{f'/fgrowth}
T(r, f'/f) = O( r \log r ) \quad \hbox{as $r \to \infty$.}
\end{equation}

If, in addition, $f$ has finite order then one of the following two conclusions holds:
\begin{eqnarray}
 \label{ReP}
&(a)& \hbox{ $f = R_1 e^{P_1}$ for some rational function $R_1$ and polynomial $P_1$;} \\
 \label{conc2}
&(b)& \quad m=2 \quad \hbox{and} \quad f(z) = \frac{A( A_1 z + 1)}{U_1(z)e^{i(B_1 z+B_2)}-1} ,
\end{eqnarray}
where $A \in \C$, while $U_1$ is a rational function with $|U_1(x)| = 1$ for all $x \in \R$,  and
$A_1,  B_1,   B_2 $ are real numbers with $B_1 \neq 0$.

Conversely, if  $f$ is as in (b) then $f$ satisfies (i), (ii) and (iii) with $m=2$. 
\end{thm}

For example, if $g(z)= z/(e^{iz}-1)$ then all but finitely many zeros of $g''$ are real by Theorem~\ref{thm1}
(see also Lemma \ref{lemrationals2}(II) below), but it is easy to check that $g'$ has infinitely many non-real zeros.
Obviously if $f$ is transcendental and is given by (a) then every derivative of $f$ has finitely many zeros. 
Examples (III), (IV) and (V)
arising from Theorem \ref{thmsnr}
show that (\ref{f'/fgrowth}) is not far from being sharp and that, at least for $m=2$, the hypothesis that $f$ has finite order
is not redundant in the second assertion of Theorem \ref{thm1}. 
Note that the analogous problem when $f$ is real was treated, but again not fully solved, in \cite{HelW3, Lawiman10,Lawiman13,rossireal}. 

The next result deals with 
strictly non-real meromorphic functions $f$ with only real zeros and poles  such that $f''/f$ is real.
Such functions  do exist but the following theorem shows that, except in one trivial case, 
the second derivative has at least one non-real zero. 

\begin{thm}
 \label{f''freal}
Let $f$ be a strictly non-real transcendental  meromorphic function in the plane, with 
finitely many zeros and poles in $\C \setminus \R$, and assume that $f''/f$ is real.
Then 
\begin{equation}
 \label{f'frepa}
\frac{f'}{f} = - \frac{\beta'}{2 \beta} + i \beta , \quad \frac{f'}{f} + \frac{g'}{g} = - \frac{\beta'}{ \beta} \, , \quad
\end{equation}
where $g = \widetilde f$ and  $\beta$ is real and meromorphic in the plane, with finitely many poles, none of them real, 
and finitely many non-real zeros.
Furthermore, $f$ has finitely many zeros. 

If, in addition, $f''$ has finitely many non-real
zeros, then $f$ satisfies~(\ref{ReP}): in particular, if all zeros and poles of $f$ and $f''$ are real 
then $f(z) = Ae^{iBz} $, where $A, B \in \C$  and $B $ is real. 
\end{thm}

 It follows from (\ref{f'frepa}) that
a zero of $\beta$ is a pole of $f$ and hence of $f''/f$,
while a pole of $\beta$ is a zero of $f$ or $\widetilde f$: 
thus if $f$ has only real zeros and $f''/f$ is entire then $\beta$ has neither zeros nor poles, and so
Theorem \ref{f''freal} contains \cite[Theorem 5]{HSW}.
Observe further that if $\beta $ is a real entire function with real zeros, all of even multiplicity,
then  (\ref{f'frepa}) defines a strictly non-real
meromorphic function $f$ with real poles and no zeros, such that $f''/f$ is real. 

\begin{cor}
 \label{cor1}
Let $H$ be a non-constant real meromorphic function in the plane with only real zeros and poles. Then any strictly non-real
meromorphic solution
in the plane of the equation $w''+Hw=0$ has at least one non-real zero.  
\end{cor}

Corollary \ref{cor1} follows at once from the last part of Theorem \ref{f''freal}, since any pole of a meromorphic solution 
of $w''+Hw=0$ is automatically a pole of $H$. 
The assertion of Corollary \ref{cor1} is not valid for real solutions, as the example $w = \tan z$, $H(z) =
-2 \sec^2 z$  immediately shows. 

The next two main results of this paper deal with the case of real functions. 
It is known  \cite{BEL,SS}  that if $f$ is a real transcendental entire function then $f$ and $f''$ have only
real zeros if and only if $f$ belongs to the Laguerre-P\'olya class $LP$, 
consisting of all entire functions which are locally uniform limits of real polynomials
with real zeros, in which case all derivatives of $f$ have only real zeros. 
For the real meromorphic case, the following was conjectured in 
\cite{HSW}. 

\begin{con}[\cite{HSW}]
 \label{conhsw}
Let $f$ be a real transcendental meromorphic function in the plane with at least one pole, and assume that all zeros and poles of
$f$, $f'$ and $f''$ are real, and that all poles of $f$ are simple. Then 
$f$ satisfies
\begin{equation}
 \label{hswform}
f(z) =  C  \tan (az+b) + Dz + E, \quad  a, b, C, D, E \in \R. 
\end{equation}
\end{con}

Without the condition that $f$ has only simple poles, there are further examples 
for which $f$, $f'$ and $f''$ have only real zeros and poles, such as $(2+\tan z)^2$ (see  \cite{toppilahsw}), as well as a substantial
collection whose existence is established  by Theorem 5 of \cite{Hin4}. 
While Conjecture \ref{conhsw} appears to be difficult to resolve in general, results  proved 
in \cite{HSW,HinRos}, and refined further in \cite{Laams09,ankara,nicks}, show in particular that the conjecture  is true subject to 
the additional hypothesis  that $f'$ omits some finite value, as is the case for the functions in (\ref{hswform}).

Theorems \ref{kreinthm} and \ref{f''=0real} below will resolve two further special cases of Conjecture \ref{conhsw}, each of them linked
to functions of the form (\ref{hswform}). 
Consider first a real  transcendental meromorphic function $f$ in the plane
which maps the open upper half-plane $H^+$ into itself: of course, $f$  also maps
the open lower half-plane $H^-$ into itself. 
Such functions $f$ have only real zeros and poles, all necessarily simple,
and by a theorem of 
Chebotarev 
\cite[Ch. VII, p.310, Theorem 2]{Le} they have a representation 
\begin{eqnarray}
f(z) &=& A z + B - \frac{d}{z} + 
\sum A_k \left( \frac1{a_k - z} - \frac1{a_k} \right), \nonumber \\
\quad 
& &
B \in \R, \quad a_k \in \R \setminus \{ 0 \}, \quad A, d, A_k \in [0, \infty) , \quad \sum \frac{A_k}{a_k^2} < \infty .
\label{krein2}
\end{eqnarray}
A well known example is  $f(z) = \tan z$. 
Conversely, any function $f$ given by  an expansion (\ref{krein2}) is real and maps $H^+$ into itself. 
This class is closely linked to the Laguerre-P\'olya class, because if $g \in LP$ then $f = -g'/g$ either is
constant or satisfies (\ref{krein2}) (see \cite{Le,SS}). 

\begin{thm}
 \label{kreinthm}
Let $f$ be a transcendental meromorphic function in the plane given by a series expansion (\ref{krein2}). 
If $m \geq 3$ then $f^{(m)}$ has infinitely many non-real zeros. 
If $f''$ has only real zeros, then $f$ satisfies (\ref{hswform}).

If $f''$ has finitely many non-real zeros, then 
\begin{equation}
 f(z) = Az+B + \frac{ R(z)e^{icz} -1}{A_1 R(z) e^{icz} - \overline{A_1}} \, , 
\label{conc1}
\end{equation}
where $A \geq 0$, $B \in \R$, $c \in (0, \infty) $, $ A_1 \in H^+$, 
and
$R$ is a rational function with all its zeros in $H^+$ and all its poles in $H^-$, and with
$|R(x)| = 1$ for all $x \in \R$.

\end{thm}

Conversely, if $f$ is given by (\ref{conc1}) with $R$ and the coefficients as in the last conclusion of Theorem \ref{kreinthm},
then $f$ maps $H^+$ into itself, and all but finitely many zeros of $f''$ are real by \cite[Lemma 3.2]{Laams09}. 
The next result in the direction of Conjecture \ref{conhsw} concerns the case where zeros of $f''$ are zeros of $f'$, 
as holds for example when $f(z) = z - \tan z$. 

\begin{thm}
 \label{f''=0real}
Let $f$ be a real transcendental meromorphic function in the plane such that:\\
(a) all but finitely many zeros and poles of $f$ and $f'$ are real;\\
(b) all but finitely many zeros of $f''$ are zeros of $f'$;\\
(c) the poles of $f$ have bounded multiplicities;\\ 
(d) either $f$ has finitely many multiple poles, or $f$ has finitely many simple poles.

Then $f$ satisfies either (\ref{ReP}) or (\ref{hswform}). 
\end{thm}

It would clearly be preferable to know whether Theorem \ref{f''=0real} holds without hypotheses (c) and (d), 
but the present method does not deliver this, and
in particular it seems difficult to exclude the possibility that $f$ has simple poles interspersed with double poles.
Of course hypothesis (d) automatically holds if $f$ is as in Conjecture \ref{conhsw}, or is
itself the derivative of a meromorphic function in the plane. 
Note that \cite[Theorem 5]{Hin4} gives rise to the example
$$
f(z) = \frac13 \tan^3 z - \tan z , \quad f'(z) = \tan^4 z - 1 , \quad f''(z) = 4 \tan^3 z \sec^2 z ,
$$
for which $f$, $f'$ and $f''$ have only real zeros and poles. Here zeros of $f''$ are zeros of $f'+1$, rather than of
$f'$, and $f$ does not satisfy
(\ref{hswform}). 
A key ingredient in the proof of Theorem \ref{f''=0real} will be to show that $f$ has finite order,
so that the following result  \cite[Theorem 3]{Lagehring} becomes relevant.

\begin{thm}[\cite{Lagehring}]\label{theorem1}
Let $f$ be a meromorphic function  in the plane  with  the following properties:
\\
(i) $f$ has finite lower order;\\
(ii) the zeros of $f'$ have bounded multiplicities;\\
(iii) all but finitely many zeros  of $f''$ are zeros of $f'$;\\
(iv) there exists  $M \in (0, + \infty) $
such that if $\zeta $ is a pole of $f$ of multiplicity $m_\zeta$ then $m_\zeta \leq M + |\zeta|^{M} $;\\
(v) there exist positive real numbers $\kappa$ and $R_0$ such that if $z$ is a zero of $f''$ 
with $|z| \geq R_0$ then  
$|f(z) - \alpha z | \geq \kappa |z|$ for all  finite non-zero
asymptotic values $\alpha$ of~$f'$. 

Then $f'' = Re^P$ with $R$ a rational function and $P$ a polynomial.
\end{thm}

Hypotheses (i) and (v) are not redundant in Theorem \ref{theorem1}, as shown by $f(z) = z - \tan z$ and examples given in \cite{Lanew}. 
The proof of Theorem \ref{f''=0real} also relies heavily on the next result. 

\begin{thm}
 \label{f''=0}
Let $n \geq 2$ be an integer, and let $f$ be a meromorphic function of finite lower order in the plane, with infinitely many poles, 
such that:\\
(i) all but finitely many zeros and poles of $f'$ have multiplicity $n$;\\
(ii) all but finitely many zeros of $f''$ are zeros of $f'$.

Then there exist $a, b, C, \lambda \in \C$ such that
\begin{equation}
 \label{f'form}
f'(z) = C \left( \frac{\lambda e^{az+b} - 1}{e^{az+b}-1} \right)^n, \quad 
aC \neq 0, \quad \lambda^n = 1, \quad \lambda \neq 1 .
\end{equation}
Furthermore, there does not exist a meromorphic function $h$ in the plane with $h'=f$.
\end{thm}

In the converse direction, it follows from Lemma \ref{lemresidue} below that 
the function in (\ref{f'form}) is indeed the derivative of a meromorphic function of finite order in the plane. 

It is worth noting that Theorem \ref{f''=0} fails completely for infinite lower order, as shown by the following example based on 
the Mittag-Leffler theorem, which is similar to Shen's construction of Bank-Laine functions with prescribed zeros \cite{Shen2}. 
Let $n \geq 2$ be an integer,
let $(a_k)$ be any complex sequence which tends to infinity without repetition, and for each $k$ let $b_k = \pm n$. 
Let $G$ be an entire function with a simple zero at each $a_k$ and no other zeros. Applying the Mittag-Leffler theorem then gives an entire function $H$
such that, for each $k$, 
$$
G(z) e^{H(z)} = \frac{z-a_k}{b_k}  + O( |z-a_k|^{n+1} ) \quad \hbox{as $z \to a_k$.}
$$
Next, a meromorphic function $g$ in the plane is determined by the formula $g/g' = G e^H$. This gives, for each $k$, as $z \to a_k$,
$$
\frac{g'(z)}{g(z)} = \frac{b_k}{z-a_k} + O( |z-a_k|^{n-1} ) , \quad 
g(z) = (z-a_k)^{b_k} (C_k +  O( |z-a_k|^{n} )) , \quad C_k \in \C \setminus \{ 0 \} .
$$
Since  $g'/g$ has no zeros by construction, the formula 
$f' = g$ now defines a meromorphic function $f$ in the plane satisfying all the hypotheses of Theorem \ref{f''=0} 
except for that of finite lower order, and each $a_k$ is a zero or pole of $f'$, depending on the sign of $b_k$. 
Moreover, $g$ is in fact the $(n-1)$'th derivative of a meromorphic function in the plane.

\section{Preliminaries}

The following theorem from \cite{FHP,La5} will be required.

\begin{thm}[\cite{FHP,La5}]
\label{thmB}
Let $f$ be a  meromorphic function in the plane, with finitely many zeros, but not of the form (\ref{ReP}). Then
$f^{(m)}$ has infinitely many zeros for every $m \geq 2$.  
\end{thm}

\begin{lem}
 \label{lem1}
Let $f$ be  a non-constant meromorphic function in the plane which satisfies at least one of the following two conditions:\\
(a)  $f$ and $f''$ have finitely many non-real zeros and poles;\\
(b) $f$ and $f^{(m)}$ have finitely many non-real zeros, for some $m \geq 3$. \\
Then the Tsuji characteristic $T_0(r, f'/f)$ in the upper half-plane satisfies
\begin{equation}
T_0(r, f'/f) = O( \log r ) \quad \hbox{as $r \to \infty$.}
\label{uhp}
\end{equation}
\end{lem}
\textit{Proof.} For details of the Tsuji characteristic see \cite{GO,Tsuji}.
Case (a)  is proved exactly as in \cite[Lemma 2.3]{BEL}, by writing
$$
F = \frac{f}{f'}, \quad F' = 1 - \frac{ff''}{(f')^2}, 
$$
so that $F$ and $F'-1$ have finitely many non-real zeros and (\ref{uhp}) follows 
from the method of Hayman's alternative \cite[Theorem 3.5, p.60]{Hay2}.
In case (b) the result is proved via  Frank's method \cite{BLa,FHP} coupled with the Tsuji characteristic.
\hfill$\Box$
\vspace{.1in}

\begin{lem}
 \label{lem2}
Let $H$ be a non-constant meromorphic function in the plane and let $G(z) = \overline{H( \bar z)}$.\\
(a)  If 
the Tsuji characteristics of $H$ and $G$ have growth given by
$$
T_0(r, H) + T_0\left(r, G \right) = O( \log r) \quad \hbox{ as $r \to \infty$,}
$$
then the Nevanlinna proximity function $m(r, H)$ satisfies 
\begin{equation}
 \label{Htsuji}
\int_R^\infty   \frac{m(r, H) }{r^3}\,  dr = O\left( \frac{ \log R }{R} \right) \quad \hbox{ as $R \to \infty$.}
\end{equation}
(b)
If $H$ satisfies (\ref{Htsuji}) and $N(r, H) = O( r \log r )$ as $r \to \infty$ then 
$T(r, H) = O( r \log r )$ as $r \to \infty$.\\
(c) If $H = e^k$, where $k$ is an entire function, and (\ref{Htsuji}) holds, then  $k$ is a  polynomial of degree~$1$. 
\end{lem}
\textit{Proof.} Applying a lemma of Levin-Ostrovskii \cite{BEL,GO,LeO} to $H$ and $G$ gives, as $R \to \infty$, 
\begin{eqnarray*}
\int_R^\infty   \frac{m(r, H) }{r^3}\,  dr  &=&
\int_R^\infty  \int_0^\pi \frac{ \log^+ |H(r e^{i \theta  } )| + \log^+ | G (r e^{i \theta  } )| }{r^3}\, d \theta \, dr \\
&\leq& \int_R^\infty \frac{T_0(r, H) + T_0(r, G)}{r^2} \, dr 
= O\left( \frac{ \log R }{R} \right) ,
\end{eqnarray*}
which proves (\ref{Htsuji}). If $H$ is as in (b) then 
(\ref{Htsuji}) holds with $m(r, H)$ replaced by $T(r, H)$ and the remaining assertions follow from the monotonicity of $T(r, H)$. 
\hfill$\Box$
\vspace{.1in}
\begin{lem}
 \label{lemrational}
Let $S$ be a rational function with $|S(x)| = 1$ for all real $x$, and let $a$ and $b$ be real numbers, with $a \neq 0$. 
Then all but finitely many solutions of $S(z)e^{i(az+b)} = 1$ are real. 
\end{lem}
\textit{Proof.} This can be deduced from \cite[Lemma 6]{nicks} but the proof is included for completeness. 
Assume that $S(\infty) = 1 =a$ and  $b=0$, and write
$g(z) = S(z)e^{iz} = e^{iz + i \phi(z)} $, in which the principal logarithm
$ \log S(z) = i \phi(z) $ tends to $0$ as $z \to \infty$, and $\phi(x) \in \R$ for real $x$ with
$|x|$ large. 
Denote by $N_{NR}$ the 
counting function of the non-real $1$-points of $g$. 
If $m \in \Z$ with $|m|$ large then the intermediate value theorem gives a solution of
the equation $x + \phi(x) = 2m \pi $ in $((2m-1) \pi, (2m+1)\pi ) $. 
Applying Nevanlinna's first fundamental theorem now yields, as $r \to \infty$,  
$$
\frac{r}\pi - O( \log r) + N_{NR}(r) \leq 
N(r, 1, g)  \leq T(r, g)+ O(1) \leq  \frac{r}{\pi} + O( \log r).
$$
\hfill$\Box$
\vspace{.1in}

\begin{lem}
 \label{lemnopoles}
Let $f$ and $g$ be meromorphic functions in the plane such that $f$, $g$ and $W = g/f$ are all non-constant.
Assume further that
\begin{equation}
 \label{identity}
\frac{f^{(m)}}{f} = \frac{g^{(m)}}{g}
\end{equation}
for some integer $m \geq 2$. If $m$ is odd then every pole of $f$ is a zero or pole of $W$. If $m$ is even then at a pole of $f$ of multiplicity $p$
which is neither a zero nor a pole of $W$, 
the function $W'$ has a zero of multiplicity $2p + m-2$, and
\begin{equation}
 \label{Nrf}
2 N_0(r, f) + (m-2) \overline{N}_0(r, f)  \leq   N(r, W/W') ,
\end{equation}
in which $N_0$ and $\overline{N}_0$ count only those poles of $f$ which are neither  zeros nor poles of $W$. 
\end{lem}
\textit{Proof.} 
Take a pole $z_0$ of $f$ of multiplicity $p$ which is neither a zero nor a pole of $W$; it may be assumed that $z_0 = 0$.
Then there exist
$\alpha$ and $\beta$ in $\C \setminus \{ 0 \}$ and a positive integer $q $ such that, as $z \to 0$, 
\begin{equation*}
 \label{laurent1}
f(z) \sim \alpha z^{-p} , \quad V(z) = W(z) - W(0) \sim \beta z^q .
\end{equation*}
The coefficient of $z^{-p+q-m}$ in the Laurent series of $U = (fV)^{(m)} - f^{(m)} V$ near $0$ is 
\begin{eqnarray*}
\alpha \beta \left[ (-p+q)\ldots (-p+q-m+1) - (-p) \ldots (-p-m+1) \right] . 
\end{eqnarray*} 
But (\ref{identity}) implies that $U $ vanishes identically, so that
\begin{equation}
 \label{preqn}
p \ldots (p+m-1) = r \ldots (r+m-1), 
\end{equation}
where $ r = p-q $. Now (\ref{preqn}) shows that $r \geq 0$ is impossible, since $r < p$, while $r < 0 \leq r+m-1$ makes the 
right-hand side of (\ref{preqn}) vanish  and so is also impossible.
Thus $s = - (r+m-1 ) > 0$, and (\ref{preqn}) yields
$$
(-1)^m s \ldots (s+m-1) = p \ldots (p+m-1) ,
$$
which forces $m$ to be even and $p = s = - (r+m-1) = - (p-q+m-1) $, so that $q = 2p + m-1$. 
\hfill$\Box$
\vspace{.1in}

Lemma \ref{lemnopoles} may be applied, in particular, if $f$ is a strictly non-real meromorphic function in the plane, with finitely
many non-real zeros and poles, such that $f^{(m)}/f $ is real for some integer $m \geq 2$: to see this, take $g(z) = \widetilde f (z)
= \overline{f( \bar z)}$.
If $m$ is odd it follows that $f$ has finitely many poles, while if $m$ is even then (\ref{Nrf}) 
yields 
$$
2 N(r, f) + (m-2) \overline{N}(r, f)  \leq   T(r, W'/W) + O(\log r ) \leq 2 m(r, f'/f) + O( \log r ),
$$
as is the case for $m=2$ and the examples mentioned following Theorem \ref{f''freal}. 

\begin{lem}
 \label{lemrationals2}
Let $T$ be a rational function with $|T(x)| = 1$ for all real $x$, let $K \not \equiv 0$ be a polynomial, 
and let $a$ and $b$ be real numbers with $a \neq 0$. 
Let 
$$
f(z) = K(z)F(z) = \frac{K(z)}{T(z)e^{i(az + b)}-1} .
$$
(I) For each $m \geq 3$ the function $f^{(m)}$ has infinitely many non-real zeros. \\
(II) If, in addition, $K$ has degree at most $1$, then all but finitely many zeros of $f''$ are real if and only if 
$L = K'/K$ is real.
\end{lem}
\textit{Proof.} It may be assumed that $a = 1$, $b =0$ and  $T(\infty) = 1$. For $|z|$ large and $\zeta \in \C$, write 
$$U(z) = iz + \log T(z), \quad f(z) = \frac{K(z)}{e^{U(z)}-1}, \quad H(\zeta) = \frac1{e^{i\zeta}-1} .$$
using the principal branch of the logarithm.

Part (I) is similar to \cite[Lemma 3.2]{Laams09}. Let $m \geq 3$, denote positive constants by $c_j$, and let 
$w$ be a non-real zero 
of $H^{(m)}$, the existence of which is assured by \cite[Lemma 3.1]{Laams09}. Take a small positive $t$ such that 
$$
 |H^{(m)}(z)| \geq c_1 \quad \hbox{and} \quad |H^{(j)}(z)| \leq c_2
$$
for $0 \leq j \leq m$ and $t \leq |z-w| \leq 3t$. Now let $n$ be a large positive integer and let 
$t \leq | z - w - 2 \pi n | \leq 3t$. Then $c_3 \leq |e^{iz}-1| \leq c_4$ and 
$$
F(z) = \frac1{e^{U(z)}-1} = \frac1{e^{iz} (1+o(1))-1} = \frac1{e^{iz} -1 + o(1)} = \frac{1+o(1)}{e^{iz} -1} = H(z) + o(1) .
$$
For $|z-w - 2 \pi n | = 2t$, applying Cauchy's estimate for derivatives yields
\begin{eqnarray*}
F^{(j)} (z) &=& H^{(j)}(z) + o(1) = O(1) \quad \hbox{for $0 \leq j \leq m$,}\\
f^{(m)}(z) &=& K(z) F^{(m)}(z) + \ldots + K^{(m)}(z) F(z) \\
&=& K(z) F^{(m)}(z) + o(1) K(z) = K(z) H^{(m)}(z) (1+ o(1)).
\end{eqnarray*}
Since $w + 2 \pi n$ is a zero of $H^{(m)}$, the assertion of part (I) now follows at once from Rouch\'e's theorem. 

To prove part (II), assume that $K'$ is constant, and write  $f'' = 2 K'F' + K F''$ and 
\begin{eqnarray}
 \label{u1}
f''(e^U-1)^3 &=& e^{2U} ( K (U')^2 - K U'' - 2K'U') + e^U ( K(U')^2 + K U'' + 2K'U') \nonumber \\
&=& e^U ( K(U')^2 + K U'' + 2K'U') (1 - Qe^U ), \nonumber \\
Q &=& \frac{     U'' + 2L U' - (U')^2}{  U'' + 2L U' + (U')^2 }, \quad L = \frac{K'}K .
\end{eqnarray}
Here $Q$ is rational but not identically zero, since $f''$ has infinitely many zeros by Theorem \ref{thmB}.
Moreover, if $x $ is real with $|x|$ large then $U'(x)$ and $U''(x)$ have zero real part, and $U'(x)^2$ is
real. If all but finitely many zeros of $f''$ are real then there exist  $x \in \R$ with $|x|$ arbitrarily large
such that $Q(x)e^{U(x)} = 1$ and so $|Q(x)| = 1$, which implies that $x$ is a zero of $Q \widetilde Q - 1$, from which it follows
that $Q \widetilde Q \equiv 1$ and $L \equiv \widetilde L$ as asserted. On the other hand, 
if  $L$ is real then $|Q(x)| = 1$ on $\R$, so that all but finitely many zeros of 
$f''$ are real by (\ref{u1}) and Lemma \ref{lemrational}. 

\hfill$\Box$
\vspace{.1in}

\begin{lem}
 \label{lemrationals}
Let $S$, $M$ and $V$ be rational functions with $S(\infty) = 1$, $M \not \equiv 0$ and $V(\infty) \neq 0$,  
and let $a$ and $b$ be complex numbers with $a \neq 0$. 
For $|z|$ large write $U(z) = az+b + \log S(z)$,
using the principal branch of the logarithm. Assume that the function $f(z)$ is meromorphic for $|z|$ large and satisfies
$$
\frac{f'}{f} = \frac{M'}{M} + \frac{V}{e^U - 1} .
$$
Then, for each $n \in \N$,
\begin{equation}
 \label{lemrat1}
\frac{f^{(n)}}f = \frac{M^{(n)}}M + \frac{V_n}{(e^U-1)^n} ,
\quad 
V_n = \sum_{j=0}^{n-1} R_{j,n} e^{jU} ,
\end{equation}
in which the coefficients $R_{j,n}$ are rational functions and satisfy, as $z \to \infty$, 
\begin{equation}
 \label{lemrat2}
R_{0,n}(z) \sim V(z)^n \quad \hbox{and} \quad R_{n-1,n}(z) \sim V(z)(-U'(z))^{n-1} .
\end{equation}
\end{lem}
\textit{Proof.} Proceeding by induction on $n$, assume that $n \in \N$ and that (\ref{lemrat1}) and (\ref{lemrat2}) both hold,
as is evidently the case for $n=1$, with $V = V_1 = R_{0,1}$. Then (\ref{lemrat1}) yields 
\begin{eqnarray*}
 \frac{f^{(n+1)}}{f} &=& \frac{M^{(n+1)}}{M} - \frac{M^{(n)}M'}{M^2} + \frac{V_n'}{(e^U-1)^n} - \frac{nV_nU'e^U}{(e^U-1)^{n+1}} \\
& & + \frac{M^{(n)}M'}{M^2} + \frac{V_n M'/M}{(e^U-1)^n} + \frac{M^{(n)}V/M}{e^U-1} + \frac{V_nV}{(e^U-1)^{n+1}} .
\end{eqnarray*}
This leads to (\ref{lemrat1}), with $n$ replaced by $n+1$ and 
$$
V_{n+1} = V_n' (e^U-1) -nV_n U' e^U + (V_n M'/M )(e^U-1) + (M^{(n)}V/M) (e^U-1)^n + V_n V .
$$
Since
$$
V_n' = \sum_{j=0}^{n-1} (R_{j,n}' + jU' R_{j,n} ) e^{jU}  ,
$$
it follows that 
$$
R_{0,n+1} = - R_{0,n}' - R_{0,n} M'/M + (-1)^n M^{(n)}V/M + R_{0,n} V 
$$
and
$$
R_{n,n+1} = R_{n-1,n}' + (n-1)U' R_{n-1,n} - nU' R_{n-1,n} + R_{n-1,n} M'/M + M^{(n)}V/M .
$$
In view of (\ref{lemrat2}) and the fact that $V(\infty) \neq 0$, this gives $R_{0,n+1}(z) \sim R_{0,n}(z) V(z) $ and
$$
 \quad R_{n,n+1}(z) = - U'(z) R_{n-1,n}(z) (1+o(1)) + o( |V(z)| ) 
\sim - U'(z) R_{n-1,n}(z) ,
$$
as $z \to \infty$, and the induction is complete. 
\hfill$\Box$
\vspace{.1in}

\begin{lem}[\cite{Lawiman12}, Lemma 4.7]\label{f'flem}
Let the function $f$ be transcendental and meromorphic in the plane and let
$k \in \N$. Let 
$E$ be an
unbounded subset of $[1, \infty )$ with the following property.
For each $r \in E$ there exist real $\theta_1(r) < \theta_2(r) \leq \theta_1(r) + 2 \pi$ 
and an arc $\Omega_r =
\{ r e^{i \theta} : \theta_1 (r) \leq \theta \leq \theta_2 (r) \}$ 
such that
\begin{equation*}
\lim_{r \to \infty, r \in E} \max \{ | z^{2k} f^{(k)}(z)/f(z) | : \, z \in \Omega_r \} = 0.
\end{equation*}
Let $N=N(r)$ satisfy
$0 \leq \log N(r) \leq o( \log r )$ as $r \to \infty$ in $ E$.
Then $f$ satisfies, for all
sufficiently large $r \in E$,
\begin{equation*}
\left| \frac{z f'(z)}{f(z)} \right| \leq kN(r)
\label{xx3} 
\end{equation*}
for all $z \in \Omega_r$ outside a union $U(r)$ of open discs having sum of radii
at most $r(k-1)/N(r)$.
\end{lem}

\begin{lem}
\label{norfamlem}
Let $k \geq 2$ and
$\rho, \sigma \in (0, \pi /2)$ and let 
$K_0 \in (0, \infty )$.
Then there exists $K_1 \in (0, \infty )$, depending only on $k$, $\rho$, $\sigma$, and $K_0$,
with the following
property. If $g$ is an analytic function on the domain 
$D = \{ z \in \C : 1/2 < |z| < 2, \, 0 < \arg z < \pi \}$ such that $g$ and $g^{(k)}$ have no zeros in $D$, 
and if 
\begin{equation*}
\min \{ |g'(e^{i\theta})/g(e^{i\theta})| : \rho \leq \theta \leq \pi - \rho \} \leq K_0,
\label{a3}
\end{equation*}
then
$|g'(e^{i\theta})/g(e^{i\theta})| \leq K_1$ for all 
$\theta \in [ \sigma, \pi - \sigma ]$.
\end{lem}

Lemma \ref{norfamlem} is standard, and follows from the fact that if $\mathcal{G}$ is the family of analytic functions on $D$
such that $g$ and $g^{(k)}$ have no zeros in $D$ then the logarithmic derivatives $g'/g$, $g \in \mathcal{G}$,
form a normal family on $D$ \cite{BLa,Schwick,Zalc}. 
The next lemma involves the
Laguerre-P\'olya class $LP$ already mentioned in the introduction~\cite{Le}.

\begin{lem}
 \label{Qlema}
Let $g \not \equiv 0$ belong to $LP$, let $M$ be a meromorphic function in the plane and write 
\begin{equation}
 \label{Qdef}
Q =  4M^3 + 6MM' + M'' , \quad Q' = 12 M^2 M' + 6 (M')^2 + 6 MM'' + M''' .
\end{equation}
(A) If  $M = R - g'/g$, in which $g$ has infinitely many zeros and $R$ is a real rational function with $R(\infty) $ finite, 
then $Q'(x)$ is positive or infinite for all $x$ in $\R$ with $|x|$ sufficiently large.\\
(B) If $M = -g'/g$ is non-constant, then 
$Q'(x)$ is positive or infinite for all $x \in \R$. 
\end{lem}
\textit{Proof.} Assume first that $M$ is as in (A) and that
$x \in \R$ with $|x|$ large. Then the standard representation \cite{Le} (see also (\ref{krein2}))
for the logarithmic derivative of a function in $LP$
leads to
$$
M' = R' - \left( \frac{g'}{g} \right)' , \quad 
M'(x) = R'(x) + 
C_0 + \sum \frac1{(a_k-x)^2} \geq \sum_{|a_k| \leq |x|} \frac1{4|x|^2} - O(|x|^{-2}) , 
$$
in which $C_0 \geq 0$ and the $a_k$ are the zeros of $g$, repeated according to multiplicity,
as well as 
$$
M'''(x) =  \sum \frac{6}{(a_k-x)^4} + O( |x|^{-4} ) \geq \sum_{|a_k| \leq |x|} \frac3{8|x|^4} - O(|x|^{-4}) .
$$
This gives 
\begin{equation}
\label{M'big}
M'(x) \sim C_0 + \sum \frac1{(a_k-x)^2} , 
\quad M'''(x) \sim   \sum \frac{6}{(a_k-x)^4}  .
\end{equation}
Write
$$
A = |M(x)|, \quad B = M'(x) > 0, \quad C = |M''(x)|, \quad D = M'''(x) > 0. 
$$
Then the Cauchy-Schwarz inequality and (\ref{M'big}) deliver
\begin{eqnarray*}
 C &\leq&  O(|x|^{-3}) + 2 \sum \left( \frac1{|a_k-x|} \cdot \frac1{|a_k-x|^2} \right)   \\
&\leq&   o ( \sqrt{ BD } ) + 2 \sqrt{ \sum \frac1{|a_k-x|^2} \sum \frac1{|a_k-x|^4} }   
\leq (1+o(1)) \sqrt{ \frac23 BD } .
\end{eqnarray*}
Assuming that the assertion of the lemma fails at $x$ gives, by (\ref{Qdef}), 
$$
12 A^2 B + 6 B^2 + D \leq 6 AC, 
$$
and squaring both sides produces
\begin{eqnarray*}
 E &=& 144 A^4 B^2 + 36B^4 + D^2 + 144 A^2 B^3 + 12 B^2 D + 24 A^2 BD \\
&\leq& 36 A^2 C^2 \leq (24 + o(1)) A^2 BD , 
\end{eqnarray*}
which implies at once that
\begin{equation}
 \label{ineq}
144 A^4 B^2 + 36B^4 + D^2 + 144 A^2 B^3 + 12 B^2 D \leq o( A^2 BD ). 
\end{equation}
But (\ref{ineq}) yields $A^4 B^2 = o( A^2 BD)$ and hence $A^2 B = o(D)$, as well as
$$0 < D^2 =  o( A^2 BD) = o(D^2),$$ 
this contradiction completing the proof of 
part (A). 

Assume now that $M$ is as in part (B) and let $x \in \R$. If $g$ has at least one zero then 
$$
M'(x) = C_0 + \sum \frac1{(a_k-x)^2} > 0 , \quad 
M'''(x) =  \sum \frac{6}{(a_k-x)^4} > 0 ,
$$
in which  $ C_0 \geq 0$ and $ a_k \in \R$,
and this time the Cauchy-Schwarz inequality gives $C^2 \leq 2BD/3$. 
If the assertion of the lemma fails at $x$ then the left-hand side of (\ref{ineq}) is non-positive,
which is impossible since $D > 0$. 

Suppose finally that $M$ is as in (B) but $g$ has no zeros. Since $M$ is assumed non-constant this forces $M' = C_0 > 0$
and $M'' = M''' = 0$,  and the conclusion of the lemma follows trivially. 
\hfill$\Box$
\vspace{.1in}

\begin{lem}
 \label{lemvalency}
Let $L$ be a real transcendental meromorphic function in the plane 
with upper half-plane 
Tsuji characteristic satisfying 
$T_0(r, L) = O( \log r )$ as $r \to \infty$, such that at least one of $L$ and $1/L$ 
has finitely many poles in $H^+$. Assume further that $F(z) = z - 1/L(z)$
has no asymptotic values $w \in H^+$, and that $F'$ has finitely many zeros in $H^+$.

Then there exists a positive integer $N$ with the following property: if $w \in H^+$ and $C$ is a component
of the set $W^+ = \{ z \in H^+: F(z) \in H^+ \}$, then each of $L$ and $F$ takes the value $w$ at most $N$ times in $C$,
counting multiplicity. 
\end{lem}
\textit{Proof.} Let $C$ be a component
of  $W^+$. The assertion concerning the valency of $F$ on $C$ is fairly standard \cite[Lemma 4.2]{BEL}:
choose a Jordan arc $\gamma^+$ which, apart from its initial point, lies in $H^+$, and is such that
every critical value $w \in H^+$ of $F$ lies on $\gamma^+$. Suppose that
$D \subseteq C$ is a component of $Y^+ = F^{-1}( H^+ \setminus \gamma^+ )$
with no non-real zero of $F'$ in $\partial D$: 
then 
the branch of
$F^{-1} $ mapping $ H^+ \setminus \gamma^+ $ to $ D$ may be analytically continued 
along $\gamma^+ \cap H^+$, giving a domain $D_1 $ with $D \subseteq D_1 \subseteq C$, mapped univalently onto 
$H^+$ by $F$, which forces $D_1 = C$. 
Thus the number
of components of $Y^+$ which lie in $C$ is bounded, independent of $C$, as is the valency of $F$ on $C$.

Controlling the number of $w$-points of $L$ in $C$, for $w \in H^+$, requires a refinement of arguments from 
\cite{Lawiman12,Lawiman13}. By \cite[Lemma 2.2]{Lawiman13}, there exist at most finitely many $\alpha \in \C$ 
such that $F(z)$ or $L(z)$ tends to $\alpha$ as $z$ tends to infinity along a path in $H^+$. This makes it possible to
choose $\theta \in (0, \pi )$ such that the two rays $P^{\pm}$, given respectively by $w = t e^{ \pm i \theta} $, $0 < t < \infty $,
contain no critical values of $L$ and no  values $\alpha $ 
such that  $L(z)$ tends to $\alpha$ as $z  \to \infty $  along a path in $H^+$.

Let $\Gamma \subseteq H^+$ be a component of $\partial C$. If $\Gamma$ is bounded, then $F$ has a pole on $\Gamma$.
On the other hand, if $\Gamma$ is unbounded, then $\Gamma$ contains a level curve 
of $F$ tending to infinity in $H^+$, on which $F(z)$ must tend to some asymptotic value belonging to $\R \cup \{ \infty \}$, because
$F$ is finite-valent on $C$. It follows that the number of components $\Gamma \subseteq H^+$ of $\partial C$ is bounded, independent of $C$.

Now take $w^* = t^* e^{i \theta} \in P^+$, and  distinct $z_1, \ldots , z_n \in C$ with $L(z_j) = w^*$. 
For each $j$, continue the branch of $L^{-1}$ mapping $w^*$ to $z_j$ along $P^+$ in the direction of decreasing $t$. 
This gives pairwise disjoint paths $\sigma_j$, which remain in $C$ since $\theta \in (0, \pi )$. Each $\sigma_j$ must tend either
to infinity or to a pole of $F$ on $\partial C$, of which only finitely many are available. Assume, after re-labelling if necessary,
that $\sigma_j$ tends to infinity for $j = 1, \ldots, m$. 

Each $\sigma_j$, for  $j = 1, \ldots, m$, may be extended to a simple path $\tau_j = \sigma_j \cup \mu_j $ in $C$, where $\mu_j$
is bounded, so that the $\tau_j$ are pairwise disjoint apart from a common starting point $z^* \in C$. After re-labelling if
necessary this gives $m-1$ pairwise disjoint domains $\Omega_j \subseteq H^+$, each bounded by $\tau_j $ and $\tau_{j+1}$.  
Because of the bound on the number of components $\Gamma \subseteq H^+$ of $\partial C$, the number of $\Omega_j$
for which $\Omega_j \not \subseteq C$  is also bounded, independent of $C$. 

Suppose now that $1 \leq k < k' \leq m-1$, and that $\Omega_k$ and $\Omega_{k'}$ are contained in $C$: then so are their closures. 
Because $F$ has no poles in $C$, the function $|L(z)|$ has a positive lower bound on the union of the $\mu_j$. Choose $q$, small and positive,
such that the circle $|w| = q$ contains no critical values of $L$ and no $\alpha$ 
such that $L(z)$ tends to $\alpha$ as $z$ tends to infinity along a path in $H^+$. 
Take $u_k \in \sigma_k$ with $L(u_k) = qe^{i \theta}$, and continue $z=L^{-1}(w)$, starting from
$qe^{i\theta}$ and along the circle $|w| = q$, so that the continuation takes $z$ into $\Omega_k$.  
Since $q$ is small, and because of the choice of $\theta$, this gives $v_k \in \Omega_k$ with $L(v_k) = qe^{- i \theta}$, 
and a simple path $\nu_k$ in $\Omega_k $ which is mapped by $L$ onto the set $\{ w = t e^{-i \theta } : 0 < t \leq q \}$. 
The fact that $L$ has no zeros in $C$ implies that $\nu_k$ must tend to infinity, and so there exists an unbounded component
$V_k$ of the set $\{ z \in \C : {\rm Im} \, (1/L(z)) > 2/q \}$, such that $V_k \cup \partial V_k \subseteq \Omega_k$. Furthermore,
the function
$$
U_k(z) = {\rm Im} \, \frac1{L(z)} \quad (z \in V_k) , \quad  U_k(z) = \frac2{q} \quad (z \not \in V_k)  ,
$$
is non-constant and subharmonic in $\C$. 
But the same argument applied to $\Omega_{k'} $ gives a corresponding component $V_{k'}$ and subharmonic function $U_{k'}$. 
A standard application of the Phragm\'en-Lindel\"of principle \cite{Hay7} yields $z$ in $V_k$ or $V_{k'}$, with $|z| $ large and 
${\rm Im} \, (1/L(z)) \geq |z|^{3/2} $, so that ${\rm Im} \, F(z) < 0$, which contradicts the fact that $z \in C$. 

Therefore  at most one of the $\Omega_j$ is contained in $C$, and this gives an upper bound, independent of $C$, for the
number $n$ of pre-images $z_j$ in $C$  of $w^* \in P^+$ under $L$. The open mapping theorem and analytic continuation of $L^{-1}$
extend this same upper bound to the number of $w$-points of $L$ in $C$, counting multiplicities, for
any $w \in H^+$. 
\hfill$\Box$
\vspace{.1in}

\begin{lem}
 \label{measurelem}
Let $Q$ be a transcendental meromorphic function in the plane such that the Nevanlinna deficiency
$\delta (\infty, Q)$ is positive. Let $C > 1$  and let $E_C \subseteq [1, \infty)$ be unbounded, such that
$T(2r, Q) \leq C T(r, Q)$ for $r \in E_C$. Let $H_r = \{ \theta \in [0, 2 \pi ] : \, 2 \log | Q( r e^{i \theta } ) | > \delta (\infty , Q) T(r, Q) \} $.
Then for large $r \in E_C$ the linear measure $m_r$ of $H_r$ satisfies
$m_r \geq d > 0$, where $d$ depends only on $C$ and $\delta (\infty, Q)$. 
\end{lem}
\textit{Proof.} This is standard. An inequality of Edrei and Fuchs 
\cite[ p.322]{EF2} yields, for large $r \in E_C$,
\begin{eqnarray*}
\frac{3\delta (\infty , Q)}4  T(r, Q)  &\leq&  m(r, Q) \leq \frac{\delta (\infty, Q)}2  T(r, Q)
+ \frac1{2 \pi} \int_{H_r} 
\log^+ | Q (r e^{i \theta })| \, d \theta \nonumber  \\
&\leq& \frac{\delta (\infty, Q)}2  T(r, Q) + 11 \left( \frac{2r}{2r-r} \right) \, m_r  \, \left( 1 + \log^+ \frac1{ m_r } \right) T(2r, Q ) \nonumber \\
&\leq&  \frac{\delta (\infty, Q)}2  T(r, Q) + 22 C  m_r  \, \left( 1 + \log^+ \frac1{ m_r } \right) T(r, Q ).
\end{eqnarray*}
\hfill$\Box$
\vspace{.1in}

\section{An auxiliary result}

The following proposition plays a fundamental role in the proof of Theorem \ref{thm1}, and in particular proves
the first assertion (\ref{f'/fgrowth}).  

\begin{prop}
 \label{keyprop}
Let the function $f$ satisfy hypotheses (i), (ii) and (iii) of Theorem \ref{thm1}.
Then 
\begin{equation}
 \label{s1}
g = \widetilde f =  R e^{ih} f = Wf , \quad \frac{g^{(m)}}{g}  = S e^{ik} \frac{f^{(m)}}{f} ,
\end{equation}
in which $\widetilde f (z) = \overline{f ( \bar z)}$, while $R$ and $S$ are rational functions, $h$ is an entire function with
\begin{equation}  
 \label{hgrowth}
T(r, h') = O( r \log r ) \quad \hbox{as $r \to \infty$,}
\end{equation}
and $k$ is a polynomial of degree at most $1$. Furthermore, $f'/f$ satisfies (\ref{f'/fgrowth}).

If, in addition, $k$ is constant in (\ref{s1}) then 
\begin{equation}
 \label{f'/fgrowth2}
N(r, f)  = O(T(r, h') + \log r ) \quad \hbox{as $r \to \infty$.}
\end{equation}
\end{prop}
\textit{Proof.} It is clear that $f$ and $g = \widetilde f$ satisfy (\ref{s1}) with $R$ and $S$ rational functions
and $h$ and $k$ entire. Now Lemma \ref{lem1} implies that, with $T_0$ the Tsuji characteristic,
\begin{equation}
T_0(r, f'/f) + T_0(r, g'/g) = O( \log r )\quad \hbox{as $r \to \infty$.}
\label{f'fg'g}
\end{equation}
Hence $h'$ and 
$e^{ik}$ satisfy the hypotheses of Lemma \ref{lem2}, from which it follows that (\ref{hgrowth}) holds, and that 
$k$ is a polynomial of degree at most $1$.

Now (\ref{f'fg'g}) also implies that (\ref{Htsuji}) holds with $H = f'/f$. But $f$ has finitely many zeros,
and so 
(\ref{f'/fgrowth}) follows provided it can be shown that  
\begin{equation}
 \label{notmanypoles}
\overline{N}(r, f) = O( r \log r ) \quad \hbox{as $r \to \infty$.}
\end{equation}
If  $k$ is non-constant then (\ref{notmanypoles}) clearly holds, since 
all but finitely many poles of $f$ are real 
$1$-points of $Se^{ik}$ by 
(\ref{s1}). 
In view of (\ref{hgrowth}), it therefore remains only to prove that (\ref{f'/fgrowth2}) holds when $k$ is constant:
if $Se^{ik} \not \equiv 1$ this follows again from (\ref{s1}). 
Suppose finally that $Se^{ik} \equiv 1$: then 
Lemma \ref{lemnopoles} may be applied, and (\ref{Nrf}) yields 
$$N(r, f) 
\leq  O( N(r, W/W')  +  \log r )  \leq O( T(r, h') + \log r )  \quad \hbox{as $r \to \infty$.} $$ 
\hfill$\Box$
\vspace{.1in}

\section{Proof of Theorem \ref{thm1}}

Let $f$ be as in the hypotheses. Since (\ref{f'/fgrowth}) has already been proved in Proposition \ref{keyprop},
it suffices to consider the case where
$f$ has finite order but (\ref{ReP}) does not hold. Then (i) and Theorem \ref{thmB} imply that $f$ has infinitely many poles and 
$f^{(m)}$ has
infinitely many zeros, all but finitely many of which are real, by (ii) and (iii). Moreover, $f$ satisfies (\ref{s1}), 
in which $R$ and $S$ are rational functions, while $h$ and $k$ are polynomials.

\begin{lem}
 \label{lemrealpoly}
It may be assumed that  $h$ and $k$ are real, and that $|R(x)| = |S(x)| = 1$ for all  $x\in \R$. 
\end{lem}
\textit{Proof.} Write $h(x) = a(x) + ib(x)$ with $a$ and $b$ real polynomials.
If $x$ is real but not a zero or pole of $f$ then $|f(x)| = |g(x)| $ and, by (\ref{s1}), 
$$
1 = |R(x)e^{ih(x)}|^2 = R(x) \overline{R(x)} \exp \left( ih(x) - i \overline{h(x)} \right) 
= |R(x)|^2 \exp( - 2 b(x) ). 
$$
Therefore $b(x) = O( \log |x| )$ 
as $|x| \to \infty$ with $x$ real. Thus $b$ is constant, and it may be assumed that $b = 0$. 
A similar argument may be applied to $Se^{ik}$. 
\hfill$\Box$
\vspace{.1in}

If $k$ is constant in (\ref{s1}) then (\ref{f'/fgrowth2}) shows that $f$ has finitely many poles, giving an immediate contradiction.
Assume henceforth that $k$ is non-constant in (\ref{s1}), and observe  that 
if $x$ is a real pole of $f$ then $S(x)e^{ik(x)} = 1$. 
Since $k$ has degree at most $1$ by Proposition \ref{keyprop}, 
it may be assumed by employing a linear change of variables that 
$S(\infty) = 1$ and $k(z) = 2 \pi z$, which, on combination with (\ref{f'/fgrowth}),
gives the following.

\begin{lem}
 \label{slem1}
The function
\begin{equation}
 \label{s2}
H(z) = \frac{f'(z)}{f(z)} \left( S(z)e^{\pi i z} - e^{- \pi i z} \right) 
\end{equation}
is meromorphic  of order at most $1$ in the plane, and has finitely many poles. 
\end{lem}
\hfill$\Box$
\vspace{.1in}


\begin{lem}
 \label{slem2}
Let $\varepsilon $ and $M$ be positive real numbers, with $\varepsilon $ small and $M$ large.
For $j=1, 2$ let $S_j$ be the sector given by $|z| \geq M$, $\varepsilon \leq (-1)^{j+1} \arg z \leq \pi - \varepsilon$. 
Then  $g$ and
$f$ satisfy, on $S_1$, 
\begin{eqnarray}
 \label{s4}
\frac{g'(z)}{g(z)} &=& T_1(z) + E(z) e^{2 \pi i z } , \quad 
\frac{f'(z)}{f(z)} 
=  W_1(z) + E(z) e^{ 2 \pi i z } , \nonumber \\
W_1(z) &=&  - \frac{R'(z)}{R(z)}  - i h'(z) + 
T_1(z) . 
\end{eqnarray}
Moreover, $f$ satisfies, on $S_2$, 
\begin{equation}
 \label{s5}
\frac{f'(z)}{f(z)} = T_2(z) + E(z) e^{-2 \pi i z } . 
\end{equation}
Here each $T_j $ is $ k_j'/k_j$ for some polynomial $k_j \not \equiv 0$ of degree at most $m-1$,  and writing
$\chi (z) = E(z)$ on $S_j$ means that $\log^+ | \chi (z)| = o(|z|) $ as $z \to \infty$ in $S_j$.
\end{lem}
\textit{Proof.} It suffices to give the proof of (\ref{s4}), that of (\ref{s5}) requiring only trivial modifications.
The function $f$ has finitely many zeros and non-real poles, and 
$k(z) = 2 \pi z$. Hence  (\ref{s1}) and standard estimates for logarithmic derivatives \cite{Gun2} show that
$$
\frac{f^{(m)}(z)}{f(z)} = E(z) , \quad g^{(m)}(z) = \delta (z) g(z) , \quad \delta (z) = E(z) e^{2 \pi i z } ,
$$
on the sector $S_1$. 
Fix a large positive $A $, and for $z \neq 0$ let  $L_z$ be 
the path  consisting of the shorter circular arc from $iA$ to 
$z^* = Az/|z|$, followed by the straight line segment from $z^*$ to $z$. If $A$ is large enough then
\begin{equation}
 \label{Adef}
2^m \int_{L_z}  |  \delta (t) t^{m-1} | \, |dt| \leq 1 
\end{equation}
for all $z \in S_1$ with $|z| \geq A$. Now there exist constants $a_j$, independent of $z$ for $z \in S_1$, such that 
$$
g(z) = a_{m-1}z^{m-1} + \ldots + a_0 + \int_{iA}^ z  \frac{(z-t)^{m-1}}{(m-1)!} \delta (t) g(t)  \, dt   ,  
$$
which can be written in the form 
$$
q(z) =  \frac{g(z)}{z^{m-1}}  = a_{m-1} + \ldots + \frac{a_0}{z^{m-1}} + \int_{iA}^z \frac{(1-t/z)^{m-1}}{(m-1)!}  \delta (t) t^{m-1} q(t) \, dt .
$$
The first step is to show that $q$ is bounded for  $z \in S_1$ with $|z| \geq A$. If this is not the case then it is possible to
choose $z \in S_1$ with $|z| \geq A$ and
$q(z) = N$ large, such that $|q(t)| \leq |N|$ on $L_z$.  
Since $ |t| \leq |z|$ on $L_z$ this gives, using (\ref{Adef}),  
$$
|N| \leq |a_{m-1}| + \ldots + |a_0| + 2^{m-1} |N|  \int_{L_z} \, |  \delta (t) t^{m-1} | \, |dt|  \leq |a_{m-1}| + \ldots + |a_0| +  \frac{|N|}2 ,
$$
which is obviously a contradiction if $N$ is large enough. 
It follows that, for $z$ in $S_1$, 
\begin{eqnarray*}
g(z) &=& k_1(z) + \int_{i \infty }^z  \frac{(z-t)^{m-1}}{(m-1)!} \delta (t) g(t)  \, dt = 
k_1(z)  + \int_{i \infty }^z  \frac{(z-t)^{m-1}}{(m-1)!} E(t) e^{2 \pi i t }   \, dt, \\
g'(z) &=&k_1'(z) + \int_{i \infty }^z  \frac{(z-t)^{m-2}}{(m-2)!} \delta (t) g(t)  \, dt , 
\end{eqnarray*}
in which the path of integration $\Lambda_z$
is  along the positive imaginary axis from $i \infty $ to $i |z|$ followed by the shorter arc of the circle $|t| = |z|$
from $i|z|$ to $z$,
while $k_1 $ is a polynomial of degree at most $m-1$. Since $|z| \leq |t|$ on $\Lambda_z$ this implies that
$$
|g(z) - k_1(z)| \leq 2^{m-1} \int_{\Lambda_z}  |t^{m-1} \delta (t) g(t) |  \, |dt| 
 \leq  \int_{\Lambda_z}  |E(t)|  e^{- 2 \pi \, {\rm Im} \, t  }   \, |dt|\quad \hbox{on $S_1$.}
$$

The next step is to show that $k_1 \not \equiv 0$. If $k_1$ vanishes identically then obviously $g(z)$ tends to~$0$ on the positive imaginary
axis, so take a large positive $y $ such that $|g(is)| \leq |g(iy)| $ for all real $s \geq y$, which gives
$$
|g(iy) \leq |g(iy)| \, 2^{m-1} \int_y^\infty  s^{m-1} |\delta(is)| \, ds ,
$$
an evident contradiction if $y$ is large enough. 

Splitting the path $\Lambda_z$ into the part from $i \infty$ to $4i |z|$ and the part $\Lambda_z^*$ from $4i |z|$ to $z$ 
now yields, for large $z$ in $S_1$, 
$$
\int_{\Lambda_z}  |E(t)|  e^{- 2 \pi \, {\rm Im} \, t  }   \, |dt| \leq 
|  e^{2 \pi i z } | \int_{\Lambda_z^*} |E(t)| \, |dt| + e^{-4 \pi |z|} \int_{4|z|}^\infty |E(is)| e^{- \pi s} \, ds 
$$
and hence 
$$
g(z) = k_1(z) + E(z) e^{2 \pi i z } , \quad 
g'(z) = k_1'(z) + E(z) e^{2 \pi i z } ,
$$
which leads to (\ref{s4}) with $T_1=k_1'/k_1$ and completes the proof of the lemma. 
\hfill$\Box$
\vspace{.1in}


It now follows from (\ref{s2}) and (\ref{s4}) that
\begin{equation}
 \label{s6}
H(z) = \frac{f'(z)}{f(z)} \left( S(z)e^{\pi i z} - e^{- \pi i z} \right)  = 
- W_1(z) e^{- \pi i z} + E(z) e^{  \pi i z } \quad \hbox{on $S_1$, }
\end{equation}
and from (\ref{s2}) and (\ref{s5}) that 
\begin{equation}
 \label{s7}
H(z) =   T_2(z) S(z) e^{ \pi i z} + E(z) e^{ - \pi i z } \quad \hbox{on $S_2$. }
\end{equation}
Since $H$ has finite order and finitely many poles, and $\varepsilon$ may be chosen
arbitrarily small,  the Phragm\'en-Lindel\"of principle gives
$$
H(z) =  T_2(z)S(z) e^{ \pi i z}  - W_1(z) e^{- \pi i z} = T_2(z) \left( S(z)e^{\pi i z} - e^{- \pi i z} \right) + V(z) e^{- \pi i z} ,
$$
in which 
$$
\quad V = T_2 - W_1 =  T_2 - T_1 + \frac{R'}{R}  + i h'  
$$
is a rational function. Using (\ref{s2}) again, this leads to 
\begin{eqnarray}
\frac{f'(z)}{f(z)} &=&  T_2(z)+  \frac{  V(z)}
{S(z)e^{2 \pi i z} - 1} .
 \label{s8}
\end{eqnarray}

Recalling that $S(\infty) = 1$ and using the principal logarithm write, for $|z| $ large, 
\begin{equation}
 \label{t1}
U(z) = 2 \pi i z + \log S(z) , \quad 
 \frac{f'}{f} = T_2 +  \frac{V}{e^U-1} .
\end{equation}

\begin{lem}
 \label{lemUV}
The  functions $U$ and $V$ satisfy $V = - U'$. 
\end{lem}
\textit{Proof.} Observe first that (\ref{t1}) shows that $f$ has infinitely many real  poles $x$ with multiplicity 
\begin{equation}
 \label{mx}
m_x = - \frac{V(x)}{U'(x)} \sim - \frac{V(x)}{2 \pi i } ,
\end{equation}
and so $V(\infty) \neq 0$. Furthermore, 
$T_2 = k_2'/k_2$, where $k_2 \not \equiv 0$ has degree at most $m-1$. 
Thus $f$ satisfies the hypotheses of Lemma \ref{lemrationals}, with $M=k_2$, by (\ref{t1}).  
It follows from (\ref{lemrat1}) and (\ref{lemrat2}) that, as $z \to \infty$ in the sector $S_2$,
on which $e^U$ is large, 
\begin{equation}
 \label{t2}
\frac{f^{(m)}(z)}{f(z)} = \frac{V_m(z)}{(e^{U(z)}-1)^m} \sim R_{m-1,m}(z)e^{-U(z)} 
\sim V(z) (-U'(z))^{m-1} e^{-U(z)} .
\end{equation}

On the other hand, since $|S| = 1$ on $\R$,
which implies that $\widetilde U' = - U'$, 
formula (\ref{s8}) leads to
\begin{equation}
 \label{t1aa}
 \frac{g'}{g} = \widetilde T_2 +  \frac{\widetilde V}{e^{-U}-1} ,
\end{equation}
in which $\widetilde T_2 = \widetilde k_2'/\widetilde k_2$. 
Since $x$ and $m_x$ are real in (\ref{mx}),
it must be the case that $\widetilde V = -V$. 
Combining Lemma \ref{lemrationals} with (\ref{t1aa}) now 
yields, as $z \to \infty$ in $S_2$,
$$
\frac{g'}{g} 
= \widetilde T_2 -  \frac{V}{e^{-U}-1}, \quad e^{U(z)} \frac{f^{(m)}(z)}{f(z)} = \frac{g^{(m)}(z)}{g(z)} 
\sim  V(z)^{m} ,
$$
using (\ref{s1}) and the fact that $e^{-U}$ is small on $S_2$. 
On comparison with (\ref{t2}) this shows that $V(\infty)/U'(\infty) $ has modulus $1$, so that $m_x$ has to be $1$
in (\ref{mx}) and the rational function $V/U' $ must be identically~$-1$.  
\hfill$\Box$
\vspace{.1in}

It now follows, using (\ref{t1}), Lemma \ref{lemUV} and the fact that $T_2 = k_2'/k_2$  for some polynomial $k_2 \not \equiv 0$, that
$f$ satisfies the hypotheses of  Lemma \ref{lemrationals2}, with 
$T(z)e^{i(az+b)} = e^{-U(z)}$ and $K/k_2$ constant. Applying Lemma \ref{lemrationals2},
part (I)  shows that $m$ must be $2$. 
Furthermore, when $m=2$, the degree of $k_2$ is at most $m-1=1$, and
part (II) of the same lemma implies that $ k_2'/k_2$ is real, so that 
any zero of $k_2$ must also be real. Conversely, if $f$ is as in conclusion (b) of the theorem, then 
all but finitely many zeros of $f''$ are real, again by Lemma \ref{lemrationals2}(II).
This completes the proof of Theorem \ref{thm1}. 
\hfill$\Box$
\vspace{.1in}

\section{Proof of Theorem \ref{f''freal}}

To prove Theorem \ref{f''freal}, assume that $f$ is a strictly non-real transcendental  meromorphic function in the plane, 
with finitely many  zeros and poles in $\C \setminus \R$, such that $f''/f$ is real. 
Write
\begin{equation}
 \label{p1}
g = \widetilde f, \quad \frac{f'}{f} = L = \alpha + i \beta , \quad \frac{g'}{g} = \widetilde L = \alpha - i \beta , \quad 
L - \widetilde L = 2i \beta , 
\end{equation}
where $\alpha$ and $\beta$ are real meromorphic functions, and $\beta \not \equiv 0$, since $g/f$ is non-constant. Then 
$$
\frac{f''}f = \alpha' + i \beta' + \alpha^2 - \beta^2 + 2 i \alpha \beta =
\frac{g''}g = \alpha' - i \beta' + \alpha^2 - \beta^2 - 2 i \alpha \beta ,
$$
from which it follows that
$$
\beta' + 2 \alpha \beta = 0, \quad L = \frac{f'}{f} = - \frac{\beta'}{2\beta} + i \beta , \quad 
\widetilde L = \frac{g'}{g} = - \frac{\beta'}{2\beta} - i \beta ,
$$
and so $f'/f$ and $\beta$ are related as in (\ref{f'frepa}). 

Now the last equation of (\ref{p1}) implies that all poles of $\beta$ are simple, and 
that $\beta$ has finitely many non-real poles. 
Moreover, 
a real pole of $\beta$ would give rise to real residues for $\beta$, $\beta'/\beta$ and $f'/f$, which is impossible by the first
equation of (\ref{f'frepa}).  
Thus $\beta$ has finitely many poles, all non-real.
It is also evident from (\ref{f'frepa}) 
that all zeros of $\beta$ have even multiplicity and are poles of~$f$, and that $\beta$ has finitely many 
non-real zeros, and finally that $f$ has finitely many zeros, as asserted. 
Obviously if $\beta$ is constant then $f(z) = Ae^{i\beta z} $, with $A$ constant.

Assume henceforth that $\beta$ is non-constant and that all but finitely many zeros of $f''$ are real. 
Then 
it is convenient to write, using (\ref{f'frepa}), 
\begin{equation}
 \label{f'frepb}
\beta = S \gamma^2, \quad P = \beta^{-1/2} ,\quad 
\frac{f'}{f} = \frac{P'}{P} + \frac{i}{P^2},  \quad M = \frac{P'}{P} =  -  \frac{S'}{2S} - \frac{\gamma'}{\gamma},
\end{equation}
where $S$ is a real rational function  and $\gamma$ is a real  
entire function  with only real zeros. 
Here $M$ is single-valued in the plane,  and $P(z)$ is single-valued for $|z|$ large,
since the zeros of $\beta$ have even multiplicity 
and the finitely many poles occur in non-real conjugate pairs.

\begin{lem}
 \label{betalem}
The function 
$\gamma$ belongs to the Laguerre-P\'olya class $LP$.
\end{lem}
\textit{Proof.} Formula (\ref{f'frepa}) and Lemma \ref{lem1} give as $r \to \infty$, using Tsuji functionals as before,
$$
m_0(r, f'/f) \leq T_0(r, f'/f) = O( \log r ), \quad T_0(r, \beta) \leq O( \log r ) + m_0(r, \beta'/\beta) 
$$
and hence $T_0(r, \beta ) = O( \log r )$, by the lemma of the logarithmic derivative for the Tsuji characteristic
\cite{GO}. Now $ \beta $ has order of growth at most $1$, by Lemma \ref{lem2}. Thus
$\gamma$ is a real entire function of order at most $1$ with 
only real zeros, and so belongs to $LP$ \cite{Le}. 
\hfill$\Box$
\vspace{.1in}

\begin{lem}
 \label{Ilem}
(a) Assume that $\gamma$ has infinitely many zeros and 
$x_0$ is a large positive real number. If $I \subseteq \R \setminus [-x_0, x_0]$ is an open interval  containing no poles of $P$, 
then $f''/f$ has at most two zeros, counting multiplicity, in $I$. \\
(b) Assume that $S = 1$ in (\ref{f'frepb}) and that $M$ is non-constant. Then $f''/f$ has at most two zeros, counting multiplicity, in
any open real interval $I$ which contains no poles of $P$. 
\end{lem}
\textit{Proof.} Observe that (\ref{f'frepb}) gives 
\begin{equation}
 \label{f''frep}
\frac{f''}{f} = \frac{P''}{P} - \frac1{P^4}  = \frac{P''}{P} - \beta^2  = \frac{P''}{P} - S^2 \gamma^4  = \frac{P^3 P'' - 1}{P^4} . 
\end{equation}
Here $P''/P$ and $P^3 P''$ are singled-valued in $\C$, since $P^2$ and $P'/P$ are.

Suppose first that $\gamma$ and $I$ are as in (a). Then $M=P'/P$ satisfies the hypotheses of part (A) of 
Lemma \ref{Qlema}, by (\ref{f'frepb}), and so the function $Q$ in (\ref{Qdef}) has at most
one zero in $I$, counting multiplicity. 
Hence the same is true of 
$$(P^3P'')' =
P^4 \left( \frac{P'''}P + 3 \frac{P'}P \frac{P''}P \right) =
P^4 \left( M^3 + 3MM' + M'' + 3 M( M^2 + M') \right) = P^4 Q. 
$$
This implies that
$P^3 P'' -1$  has at most two zeros in $I$, counting multiplicity, and so has $f''/f$, by (\ref{f''frep}). 
Part (b) is proved the same way, since if $S=1$ and $M$ is non-constant then $M$  satisfies the hypotheses of 
Lemma \ref{Qlema}(B).
\hfill$\Box$
\vspace{.1in}

\begin{lem}
 \label{betapoly}
The function $\beta$ is rational, and $f$ satisfies (\ref{ReP}). 
\end{lem}
\textit{Proof.} Assume that  $\beta$ is  transcendental. If $\beta$ has finitely many zeros then 
$\beta (z) = R_1 (z) e^{b_1z} $, with $R_1$ a rational function and $b_1 \in \R \setminus \{ 0 \}$, 
and (\ref{f''frep}) shows that $f''/f$ has infinitely many non-real zeros,
which is a contradiction.

Assume henceforth that $\beta$ has infinitely many zeros; then so has $\gamma$. 
Since $f''/f$ has a double pole at each real pole $x$ of $P$ with $|x|$ large, and has
finitely many non-real zeros, Lemma \ref{Ilem}(a)
implies that the following estimates hold as $r \to \infty$. First, 
$$
n(r, f/f'') \leq n(r, f''/f) + O(1), \quad N(r, f/f'') \leq N(r, f''/f) + O(\log r) ,
$$
from which applying Jensen's formula yields, in view of (\ref{f'frepb}), (\ref{f''frep}) and the fact that $\beta$ has finite order,
\begin{eqnarray*}
 2 m(r, \beta ) &\leq& m(r, f''/f) + O( \log r ) \leq m(r, f/f'') + O( \log r )  \\
&\leq& T(r, f''/f) + O( \log r )  = O( T(r, \beta ) ) .
\end{eqnarray*}
Thus the zeros of $f''/f$ have positive Nevanlinna deficiency  $\delta (0, f''/f)$. 

A contradiction will now be obtained using a method similar to the proof of \cite[Lemma 5.4]{Lawiman12}.
Since $\beta$ and $f''/f$  have finite order, 
a well known result of Hayman \cite[Lemma 4]{Hay3}
gives $C_1 > 0$ and a set $E_1 \subseteq [1, \infty)$, of
positive lower logarithmic density, such that  
\begin{equation}
\label{c4}
T(4s, \beta) \leq  C_1 T(s, \beta ) \quad \hbox{and}  \quad T(4s, f''/f) \leq C_1 T(s, f''/f)   
\end{equation}
for $s \in E_1$. By estimates from \cite{Gun2}, the function $\beta$ also satisfies
\begin{equation}
 \label{c4a}
\left| \frac{\beta'(z)}{\beta(z)} \right| \leq r^{M_0} \quad \hbox{for $|z|=r \not \in F_2$},
\end{equation}
where $M_0$ is a positive constant and $F_2$ has finite logarithmic measure. 

Now let $\sigma$, $K_0$, $K_1$ and $K_2$ be positive
constants, with $K_0$, $K_1/K_0$ and $K_2/K_1$ large, and $\sigma$ small.
Let $s \in E_1$ be large. Since $f''/f$ is transcendental and $\delta (0, f''/f) > 0$, a standard application of (\ref{c4})
and Fuchs' small arcs lemma \cite[p.721]{Hay7} 
gives $r \in [s, 2s] \setminus F_2$ and 
an arc of the circle $|z| = r$, of angular measure $6 \sigma$, on which
$|f''(z)/f(z)| \leq r^{-5} $. The fact that $f''/f$ is real then implies that 
$|f''(z)/f(z)| \leq r^{-5} $ on a subarc $I_r$ of 
$\{ z \in \C : |z| = r, \sigma \leq \arg z \leq \pi - \sigma \}$ of angular measure at least $\sigma$. 
Next, applying Lemma \ref{f'flem} with $k=2$ and $N(r)=K_0$ 
shows that there exists $z \in I_r$ with $|zf'(z)/f(z)| \leq K_1$.
Now Lemma \ref{norfamlem}, applied to the function $f(rz)$, delivers
$|zf'(z)/f(z)| \leq K_2$ for all $z$ with $|z| = r$, $ \sigma \leq \arg z \leq \pi - \sigma $.
Because $\beta$ is real, combining this estimate with (\ref{f'frepa}), (\ref{c4}) and
(\ref{c4a}) yields an unbounded set of positive $r$ such that 
$T(2r, \beta ) \leq T(4s, \beta) \leq C_1 T(s, \beta) \leq C_1 T(r, \beta)$, 
and such that 
$| \beta (z)| \leq r^{M_0}$ for all $z$ with $|z| = r$, apart from  a set $J_r$ of
angular measure at most $4 \sigma$, where  $\sigma$ may be chosen arbitrarily small, independent of $C_1$. 
Since $\beta$ has finitely many poles, this contradicts Lemma \ref{measurelem}.

Thus  $\beta$ is rational, as asserted, and so is $f'/f$ by (\ref{f'frepa}), which implies (\ref{ReP})  and completes
the proof  of the lemma. 
\hfill$\Box$
\vspace{.1in}

To finish the proof of the theorem assume henceforth that all zeros and poles of
$f$ and $f''$ are real.
Then $\beta$ has no poles, by (\ref{p1}), and so it may be assumed that $S = 1$ in (\ref{f'frepb}). 
Since zeros of $\beta$ have even multiplicity,
and the case where $\beta$ is constant has already been disposed of,
it can  now be assumed  that $\beta$ is a polynomial with real zeros, of even positive degree, and $M$ is
non-constant in (\ref{f'frepb}). 
Thus (\ref{f'frepa}) and (\ref{f''frep}) show that
$f''/f$ is a rational function with double poles at the zeros of $\beta$,
which are real poles of $P$. Moreover $f''/f$ has only real zeros, and by Lemma \ref{Ilem}(b) the number of zeros of $f''/f$ exceeds the
number of poles by at most $2$. Hence $f''/f$ has at most a double pole at infinity and so $\beta$ has degree at most $1$, 
by (\ref{f''frep}) again, which
is a contradiction. 
\hfill$\Box$
\vspace{.1in}

\section{Some applications of harmonic measure}

\begin{lem}[\cite{EFacta,Nev}]\label{EDREI-FUCHS}
 Let $G$ be a domain bounded by a Jordan curve $C$ consisting of a Jordan arc $B$ and
its complement $A = C \setminus B$. Let $L$ be a rectifiable curve in $G$ joining $a \in A$
to $b \in B$, and for $z \in L$ let $\rho(z)$ be the distance from $z$ to $A$. 
Then the harmonic measure $\omega(z)$ of $B$ with respect to $G$ satisfies, for $z $ on $ L$,
$$
\omega(z) \geq \frac1{2\pi} \exp \left( -4 \int_z^b \, \frac{|du|}{\rho(u)} \right) ,
$$
in which the integration is from $z$ to $b$ along $L$. 
\end{lem}

\begin{lem}
 \label{lemhfn}
Let $Q$ be a transcendental meromorphic function of finite order in the plane, such that the zeros of $Q$ have positive
Nevanlinna deficiency $\delta (0, Q)$. 
Assume that for each $\delta > 0$ there exists $ N(\delta) > 0$ such that 
\begin{equation}
 \label{h1}
\log |Q(z)| \leq  N(\delta) \log |z| 
\end{equation}
for all $z$ with $|z|$ large and $\delta \leq | \arg z | \leq \pi - \delta $. 

Let $\eta$ and $\varepsilon$ be positive. Then, for all sufficiently large $r$, the function $Q$ satisfies 
\begin{equation}
 \label{h2}
\log |Q(z)| \leq  2 N( \varepsilon /2) \log r  - r^{-\eta}  T(r, Q)  
\end{equation}
for all $z$ in at least one of the arcs 
$$
I^+(r, \varepsilon) = \{ r e^{i \theta} : \varepsilon \leq \theta \leq \pi - \varepsilon \}, \quad 
I^-(r, \varepsilon) = \{ r e^{- i \theta} : \varepsilon \leq \theta \leq \pi - \varepsilon \}.
$$
\end{lem}
\textit{Proof.} The initial steps are standard. 
Choose $\delta > 0$, small compared to $\eta$. 
By the same result of Hayman \cite[Lemma 4]{Hay3} as used in the proof of Lemma \ref{betapoly},
there exists $C_1 > 0$, depending on $\delta$ and the order of $Q$,
as well as a set $E_\delta \subseteq [1, \infty)$, of lower logarithmic density at least
$1-\delta/2$, such that if $s \in E_\delta$ then $T(4s, Q) \leq C_1 T(s, Q)$. 
Let $s \in E_\delta$ be large, let $H_s = \{ \theta \in [0, 2 \pi ] : \, 2 \log | Q( 2s e^{i \theta } ) | < - \delta (0, Q) T(2s, Q) \} $,
and let $m_s$ be the linear measure of $H_s$. 
Then Lemma \ref{measurelem} yields
$m_s \geq 16 \delta_1 > 0$, where 
$\delta_1 $ is small but independent of $s$. 

Now let $r$ be large and positive: then there exists $s \in E_\delta$ with
\begin{equation}
 \label{sdef}
2r \leq  s \leq   r^{1+\delta} \leq r^2 .
\end{equation}
Since $H_s$ has measure $m_s \geq 16 \delta_1$, it may be assumed without loss of generality that $Q$ satisfies
$2 \log | Q( z ) | < - \delta (0, Q) T(2s, Q)$ for all $z$ in a subset $I_s$ of $I^+ (2s, 2 \delta_1 ) $, of
angular measure at least $4 \delta_1$. Let $D_s$ be the domain
$$
\{ z \in \C : \,  s/2 < |z| < 2s, \, \delta_1 < \arg z < \pi - \delta_1 \} 
$$
and let $w \in I^+(s, \pi/4)$. Then the harmonic measure $\omega(w, I_s, D_s)$ of $I_s$ with respect to $D_s$ 
is bounded below by a positive constant $\delta_2$ which is independent of $s$ and $r$. Thus (\ref{h1}) and
the two constants theorem \cite{Nev} yield, since $Q$ is transcendental and $r$ and $s$ are large, 
\begin{equation}
 \label{h2a}
\log |Q(w)| \leq N( \delta_1) \log 2s - \frac{\delta_2 \delta (0, Q) }2 T(2s, Q)
\leq - \frac{\delta_2 \delta (0, Q) }4 T(2s, Q)  
\end{equation}
for all $w \in I^+(s, \pi/4)$. 
Next, let 
$\Omega$ be the domain 
$$
\{ z \in \C : r/2 < |z| < s, \, \varepsilon /2  < \arg z < \pi - \varepsilon /2  \} ,
$$
and let $z_0 \in I^+(r, \varepsilon )$. 
Join $z_0$ to $is$ by the simple path $\gamma$ consisting of the shorter arc of the circle 
$|z| = r$ from $z_0$ to $ir $,
followed by the radial segment $z = ix  $, $r \leq x \leq s$. Let $B = I^+(s, \pi/4)$ and
$A = \partial \Omega \setminus B$. 
Denoting by $\rho(u)$  the distance from $u$ to $A$ then gives, on integrating with respect to arc length and using (\ref{sdef}),  
$$
\int_\gamma \, \frac{|du|}{\rho(u)} \leq d_1 \left( \frac1\varepsilon + \int_r^s \, \frac{dt}{t} \right) 
\leq d_1 \left( \frac1\varepsilon + \delta \log r  \right) ,
$$
where $d_1 > 0$ is independent of $\varepsilon$, $\delta $ and $r$. This time the two constants theorem  delivers,
in view of (\ref{h1}), (\ref{h2a}) and Lemma \ref{EDREI-FUCHS}, 
$$
\log |Q(z_0)| \leq 2 N( \varepsilon /2) \log r 
-   \frac{\delta_2 \delta (0, Q) }{8 \pi}  T(2s, Q)  \exp \left( -4 d_1 \left( \frac1\varepsilon + \delta \log r  \right) \right).
$$
Since $r$ is large and $\delta/\eta$ is small, (\ref{h2}) follows for $z = z_0$, and the proof is complete. 
\hfill$\Box$
\vspace{.1in}

\begin{lem}
 \label{lemcarleman}
Let $u$ be a non-constant continuous subharmonic function in the plane, of finite order $\rho$, and let $\varepsilon > 0$. 
Let $F$ be the set of $r \in [1, \infty)$ for which there exists an arc of the circle $|z| = r$, of length at 
least $\varepsilon r$, on which $u(z) > 0$. 
Then $F$ has lower logarithmic density at least $1-  \varepsilon \rho/\pi$.
\end{lem}
\textit{Proof.} This is a standard application of a well known estimate for harmonic measure
\cite{Tsuji1}. For $r > 0$ let
$B(r, u) = \max \{ u(z): |z| = r \}$ and let $r \theta(r)$ be the length of the longest open arc of the circle $|z| = r$ on which $u(z) > 0$,
except that $\theta(r) = \infty$ if $u(z) > 0$ on the whole circle. Then, as $r \to \infty$, by \cite[p.116]{Tsuji1},
$$
\int_{[1, r] \setminus F} \, \frac{dt}{t} \leq  
\frac{\varepsilon}{\pi}  \int_1^r \, \frac{\pi dt}{t \theta (t)} \leq \frac{\varepsilon}{\pi} \log B(2r, u) + O(1) \leq \frac{\varepsilon}{\pi} 
(\rho + o(1)) \log r .
$$
\hfill$\Box$
\vspace{.1in}

\begin{lem}
 \label{nodirect}
Let $G$ be a transcendental meromorphic function of finite order  in the plane, and assume that there exist 
$\alpha_1, \alpha_2 \in \C$, not necessarily distinct, with the following property: for each $\varepsilon > 0$, the function $G$ satisfies
$G(z) \to \alpha_j$ as $z \to \infty$ with $\varepsilon < (-1)^j \arg z < \pi - \varepsilon $.
If $\beta \in ( \C \cup \{ \infty \} ) \setminus  \{ \alpha_1, \alpha_2 \}$ then the inverse function of $G$ cannot have 
a direct transcendental singularity over $\beta$.
\end{lem}
\textit{Proof.} This is again standard: for the terminology see \cite{BE,Nev}. 
Assuming without loss of generality that $G^{-1}$ has a direct transcendental singularity over
$\beta \in  \C \setminus  \{ \alpha_1, \alpha_2 \}$ gives a small positive
$\delta$, a component 
$U$ of the set $\{ z \in \C : | G(z) - \beta | < \delta \}$, and a 
continuous, subharmonic, non-constant function $u$ of finite order in the plane which satisfies
$u(z) = \log (\delta/|G(z)- \beta|)$ on $U$ and vanishes outside $U$.
Here $\delta$ may be chosen arbitrarily small, as may $\varepsilon $. But then the intersection of $U$ with the set
$\{ z \in \C : \varepsilon < | \arg z | < \pi - \varepsilon \}$ is bounded, which contradicts Lemma \ref{lemcarleman}. 
\hfill$\Box$
\vspace{.1in}

\section{Proof of Theorem \ref{kreinthm}}
Let $f$ be a transcendental meromorphic function given by (\ref{krein2}). 

\begin{lem}
 \label{lemk1}
Let $n$ be a non-negative integer and let $N_R(r, 1/f^{(n)})$ count the real zeros of $f^{(n)}$, with respect to multiplicity.
If $n$ is odd then $N_R(r, 1/f^{(n)}) = 0$. If $n$ is even then $f^{(n)}$ has at most one zero in any open interval of the real axis 
which contains no poles of $f$, and 
$N_R(r, 1/f^{(n)} ) \leq N(r, f) + O( \log r )$ as $r \to \infty$. 
Furthermore, if $a_k$ and $a_{k+1}$ are poles of $f$, with $a_k < a_{k+1}$ and no poles of $f$
in $I_k = (a_k, a_{k+1})$, then $I_k$ contains precisely one zero of $f''$. Finally, 
$m(r, f) = O( \log r)$ as $r \to \infty$. 
\end{lem}
\textit{Proof.}
The first three assertions follow from differentiating 
(\ref{krein2}), which shows that if $m$ is an odd positive integer then $f^{(m)}(x)$ is positive or infinite for every real $x$.
Next, the fact that all residues of $f$ are negative, while all poles of
$f''$ have multiplicity $3$, forces $f''$ to change sign on $I_k$. Hence $f''$ has precisely one zero in $I_k$, since $f'''$ has none there. 
The bound on $m(r, f)$ holds since $f$ is real and maps the upper half-plane $H^+$ into itself, so that 
\cite[Ch. I.6, Thm. $8'$]{Le}
\begin{equation}
\frac{1}{5}|f(i)|\frac{\sin\theta}{r}<|f(re^{i\theta}
)|<
5|f(i)|\frac{r}{\sin\theta} \quad \hbox{for}  \quad r \geq 1,\,
\theta\in (0,\pi).
\label{C1}
\end{equation}
\hfill$\Box$
\vspace{.1in}

\begin{lem}
 \label{lemk2}
Let $m \geq 3$, let $\varepsilon$ be small and positive and let $N_{NR}(r, 1/f^{m)})$ count the non-real zeros of $f^{(m)}$.
Then $f$ satisfies $(m-2 - \varepsilon ) T(r, f) \leq N_{NR}(r, 1/f^{(m)})$ as $r \to \infty$ outside a set of finite measure.
In particular, $f^{(m)}$ has infinitely many non-real zeros. 
\end{lem}
\textit{Proof.} 
Since  $f$ is transcendental with only real poles, all of which are simple,
Lemma \ref{lemk1} and 
an inequality of Frank, Steinmetz and
Weissenborn \cite{Frankgraz} (see also \cite{FW,FW2,Steiw}) yield, for large $r$ outside a set of finite measure, 
\begin{eqnarray*}
 (m+1) T(r, f) &=& (m+1) N(r, f) + O( \log r) 
= N(r, f^{(m)}) + o( T(r, f)) \\
&\leq&  N(r, 1/ f^{(m)}) + (2+\varepsilon/2) N(r, f) + o( T(r, f)) \\
&\leq& N_{NR}(r, 1/f^{(m)}) + (3+\varepsilon/2) N(r, f) + o( T(r, f)) .
\end{eqnarray*}
\hfill$\Box$
\vspace{.1in}

Lemma \ref{lemk2}  proves the first assertion of Theorem \ref{kreinthm}. 
Assume henceforth that $f''$ has finitely many non-real zeros. 
Clearly all zeros of $f'$ are non-real by Lemma~\ref{lemk1}. Let 
\begin{equation}
 \label{Fdef}
F(z) = z - \frac{f(z)}{f'(z)} , \quad 
W^+ = \{ z \in H^+: F(z) \in H^+ \}, \quad
W^- = \{ z \in H^+: F(z) \in H^- \}.
\end{equation}
It may be assumed that $A= B=0$ in (\ref{krein2}), since $f(z) - Az-B$ has the same second derivative as $f$. 
\begin{lem}
 \label{lem4a}
Let $\varepsilon > 0$. Then $f(z)/z \to 0$ as $z \to \infty$ with 
$\varepsilon < | \arg z | < \pi - \varepsilon$. 
\end{lem}
\textit{Proof.} This is standard. Fix $\delta > 0$ and let  $R \geq 1$. 
Then (\ref{krein2}) gives a rational function $T_R$, with $T_R(\infty) = 0$, such that,  
for $\varepsilon < | \arg z | < \pi - \varepsilon$, 
$$
\frac{f(z)}{z} = T_R(z) + \sum_{|a_k| > R} \frac{A_k}{a_k(a_k-z)} , \quad \left| \frac{f(z)}{z} \right| 
\leq |T_R(z)| + \sum_{|a_k| > R} \frac{A_k}{a_k^2 \sin \varepsilon } = |T_R(z)| + S.
$$
Now choose $R$ so large that (\ref{krein2}) gives $S < \delta$, and  $|z|$ so large that $|T_R(z)| < \delta$. 
\hfill$\Box$
\vspace{.1in}

\begin{lem}
 \label{lem3}
All poles of $F$ are  non-real, 
while all but finitely many zeros of $F'$ are real. In any open interval of the real axis which contains
no poles of $f$, the function $F'$ has at most two zeros, counting multiplicity. 
\end{lem}
\textit{Proof.} These assertions all follow from  Lemma \ref{lemk1} and the formula $F' = (ff'')/(f')^2$. 
\hfill$\Box$
\vspace{.1in}
\begin{lem}
 \label{lem4}
The Tsuji characteristic of $f'/f$ satisfies (\ref{uhp}), and $f$ has order of growth at most~$1$ in the plane. 
\end{lem}
\textit{Proof.} The first assertion follows from Lemma \ref{lem1}. Alternatively, it may be observed that the function
$(f-i)/(f+i)$ has modulus less than $1$ on $H^+$.

To prove that $f$ has order at most $1$, 
the function $f''/f$ will be written as a quotient as follows. 
Assume that the $a_k$ in (\ref{krein2}) are ordered so that $a_{k} < a_{k+1}$ for each $k$. 
If $|k| \geq k_0$, 
where $k_0$ is large, then $a_k$ and $a_{k+1}$ have the same sign, and by Lemma \ref{lemk1} there is precisely one zero $b_k$ of $f''$ 
in $(a_k, a_{k+1})$, counting multiplicity. Write  
$$
\psi (z) = \prod_{|k| \geq k_0} \frac{1-z/b_k}{1-z/a_k}, \quad 
0 < \sum_{|k| \geq k_0} \arg 
\frac{1-z/b_k}{1-z/a_k} =
\sum_{|k| \geq k_0} \arg 
\frac{b_k -z}{a_k -z} < \pi \quad \hbox{for} \quad z \in H^+ .
$$ 
The product $\psi$  converges by the alternating series test, and  $\psi (H^+) \subseteq H^+$.
Next, write $f''/f = \psi/g$, where $g = \psi f/f''$ has finitely many poles, using Lemma \ref{lemk1},
and all but finitely many poles of $f$ are simple zeros of $g$. 
It follows from (\ref{uhp}) and standard properties of the Tsuji characteristic that the hypotheses of Lemma \ref{lem2}(a) are satisfied
with $H = f/f''$ (and so $\widetilde H = H$). This gives (\ref{Htsuji}) with $H = f/f''$.

Now  $m(r, f) = O( \log r)$ by Lemma \ref{lemk1}, and the same is true with $f$ replaced by $\psi$, because $\psi (H^+) \subseteq H^+$. 
Therefore (\ref{Htsuji}) also holds with $H = g$. Thus Lemma \ref{lem2}(b) shows that $T(r,g)$ has order of growth at most $1$, and hence
so have $N(r, f)$ and $T(r, f)$. 
\hfill$\Box$
\vspace{.1in}

\begin{lem}
 \label{lem6a}
There does not exist $\beta \in \C \setminus \{ 0 \}$ such that $f(z)/z \to \beta$ as $z$ tends to infinity on a path in $\C \setminus \R$.  
\end{lem}
\textit{Proof.} If such an asymptotic value $\beta$  exists then the inverse function of $f(z)/z$ has a direct transcendental singularity 
over $\infty$, by Lemma \ref{lem4a}. 
But this is impossible, by Lemmas \ref{nodirect} and \ref{lem4a} and the fact that $f$ 
has finite order of growth.
\hfill$\Box$
\vspace{.1in}

\begin{lem}
 \label{lem5}
Let $\alpha \in \C \setminus \R$. 
Then the inverse function $F^{-1}$ has no direct transcendental singularities over $\alpha$. 
\end{lem}
\textit{Proof.} 
Assume that $F^{-1}$ does have a direct transcendental singularity over $\alpha \in \C \setminus \R$. 
Then, without loss of generality, there exist $\delta > 0$ and a component 
$U \subseteq H^+$ of the set $\{ z \in \C : |F(z) - \alpha | < \delta \}$, such that the function
\begin{equation}
 \label{subhdef}
u(z) = \log \frac\delta{|F(z)-\alpha|} \quad (z \in U), \quad u(z) = 0 \quad (z \in \C \setminus U),
\end{equation}
is subharmonic and non-constant in the plane. By a result of Lewis, Rossi and Weitsman \cite{LRW} there exists a path $\Gamma$ tending
to infinity in $U$ on which $u(z) \to + \infty$ with
\begin{equation}
 \label{lrw}
\int_\Gamma e^{-u(z)} \, |dz| < \infty .
\end{equation}
For $z \in \Gamma$ with $|z|$ large write
$$
z - \frac{f(z)}{f'(z)} = F(z) = \alpha + p(z) ,
\quad \frac{f'(z)}{f(z)} = \frac1{z - \alpha } + q(z), \quad |q(z)| \leq |p(z)|  = \delta e^{-u(z)}  .
$$
Hence (\ref{lrw}) shows that 
there exists a non-zero complex number $\beta$ such that $f(z) \sim \beta ( z - \alpha )$ as $z \to \infty$ 
on~$\Gamma$,  contradicting Lemma \ref{lem6a}. 
\hfill$\Box$
\vspace{.1in}

\begin{lem}
 \label{lem6}
The function $F$ has finitely many critical values, and no asymptotic values, in $\C \setminus \R$.

\end{lem}
\textit{Proof.} The fact that  all but finitely many
critical values of $F$ are real is an immediate consequence of Lemma \ref{lem3}. 
Since all  poles of $f'/f$ are real, it follows from Lemma \ref{lem4} and
\cite[Lemma 2.2]{Lawiman13}
that $F$ has finitely many asymptotic values in 
$\C \setminus \R$. Because $F$ has finite
order, 
any non-real finite asymptotic value of $F$ must give rise to a direct singularity
of $F^{-1}$, by \cite{BE}, contradicting Lemma \ref{lem5}. 
\hfill$\Box$
\vspace{.1in}

\begin{lem}
 \label{componentlem}
There exists a positive integer $M$ such that if $C$ is a component of $W^+$ or $W^-$ then $F$ takes each value at most $M$ times
in $C$, counting multiplicity. Furthermore, a component of $W^+$ (respectively, $W^-$) which contains no zeros of $f''$
is simply connected
and conformally equivalent to $H^+$ (respectively, $H^-$) under $F$, and this is true for all but finitely many components
of $W^+$ (respectively, $W^-$). 
\end{lem}
\textit{Proof.} The first assertion is proved as in Lemma \ref{lemvalency}, using Lemma \ref{lem3}; the rest is standard.
\hfill$\Box$
\vspace{.1in}

\begin{lem}
 \label{lem8}
Let $C$ be a component of $W^+$ or $W^-$ which contains no zeros of $f''$, and let $\alpha \in \R$. Then there exists  $z$ in the finite boundary
$ \partial C$ with $F(z) = \alpha$.  
\end{lem}
\textit{Proof.} Let $C$ and $\alpha$ be as in the hypotheses and assume that $\alpha \not \in F ( \partial C)$. 
Let 
$G(z) = 1/(\alpha - F(z))$, so that $G$ is univalent on $C$, and $G(C)$ is $H^+$ or $H^-$.   Let 
$g: G(C) \to C$ be the inverse function of $G$, and let $\Gamma$ be the path in $ G(C)$ given by
$$
w = it , \quad t \in \R, \quad 1 \leq |t| < \infty . 
$$
Then $\gamma=g(\Gamma)$ is a curve in $C$ on which $iG$ is real, and $\gamma$ tends either to infinity or to an $\alpha$-point of $F$ on $\partial C$.
Hence $\gamma$ must tend to infinity in $C$. For $z \in \gamma$ with $|z| $ large write
$$
z - \frac{f(z)}{f'(z)} = F(z) = \alpha - \frac1{ G(z)} = \alpha + o(1) ,
$$
which leads to
$$
\frac{f'(z)}{f(z)} = \frac1{z - \alpha +1/G(z)} = \frac1{z-\alpha} + h(z) ,
\quad \hbox{where} \quad h(z) =  O \left( \frac1{|z|^2 |G(z)| } \right)  .
$$
But Koebe's $1/4$ theorem applied to $\log g$  gives $g'(w)/g(w) = O( 1/|w|)$ on $\Gamma$ and so 
$$
\int_\gamma |h(z)| \, |dz| = \int_\Gamma O \left( \frac{|g'(w)|}{|g(w)|^2 |w| } \right)  \, |dw| =
\int_\Gamma O \left( \frac{1}{|w|^2 |g(w)| } \right)  \, |dw| < \infty. 
$$
It follows that there exists a non-zero complex number $\beta$ such that $f(z) \sim \beta ( z - \alpha )$ as $z \to \infty$ 
on~$\gamma$, and this contradicts Lemma \ref{lem6a}. 
\hfill$\Box$
\vspace{.1in}

\begin{lem}
 \label{lem11a}
Let $a \in \R$ be a zero of $f''$. Then $f$ has at least one pole in each of $(-\infty, a)$ and $(a, \infty)$. 
\end{lem}
\textit{Proof.} Suppose that $f$ has no poles in $(-\infty, a)$. Then $(X - a)^3 > 0$ 
for  every pole $X$ of $f$, and the  series expansion for
$f''$ obtained from (\ref{krein2}) shows that $a$ cannot be a zero of $f''$. 
\hfill$\Box$
\vspace{.1in}

\begin{lem}
 \label{lem11b}
Every pole of $f$ lies on the boundary of a component of $W^+$ but not in the closure of $W^-$. 
\end{lem}
\textit{Proof.} This holds because every pole $X$ of $f$ is a real fixpoint of $F$ with $F'(X) > 1$. 
\hfill$\Box$
\vspace{.1in}

\begin{lem}
 \label{lem11cc}
Let $a \in \R$ be a multiple zero of $F'$. Then $F'''(a) > 0$.
\end{lem}
\textit{Proof.} 
Lemma \ref{lemk1} shows that $a$ must be a common zero of $f$ and $f''$, and a triple zero of $F - F(a)$. 
Assume that $F'''(a)$ is negative and let $\delta$ be small and positive: 
then $a - \delta $ and $a + \delta $ both lie in $\partial W^-$. 
Let $A$ and $B$ be the nearest poles of $f$ to $a$ in $(-\infty, a)$ and
$(a, \infty)$ respectively; these exist by Lemma \ref{lem11a}, and Lemma \ref{lem11b} ensures that
each lies on the boundary of a component of $W^+$. It follows that $F$ must have critical points in $(A, a)$ and $(a, B)$,
contradicting Lemma~\ref{lemk1}. 
\hfill$\Box$
\vspace{.1in}

\begin{lem}
 \label{lem11}
The function $f'$ has finitely many zeros, and none at all 
if $f''$ has only real zeros. 
\end{lem}
\textit{Proof.} Let $w$ be a zero of $f'$. Then $w$ is non-real by Lemma \ref{lemk1},
and it may be assumed that $w \in H^+$. Thus $w$ is a pole of $F$: with finitely many exceptions, and none at all if $f''$ has
only real zeros, the pole of $F$ at $w$ is simple. 

Assume henceforth that $w \in H^+$ is a zero of $f'$ and a simple pole of $F$: then $w$ lies on the boundary of
a uniquely determined component $C_w$ of $W^-$. 
Consider those $w$ for which the component $C_w$ either is multiply connected, 
or has a non-real zero of $f''$ in its closure. There are finitely many of these, by Lemma \ref{componentlem},
and none if $f''$ has only real zeros. 

Attention may thus  be restricted to those $w$ for which
$C = C_w$  is simply connected, with no non-real zero of $f''$ in 
its closure. Then $F$ maps $C$ univalently onto $H^-$, and
$F(\partial C) = \R \cup \{ \infty \}$, by Lemmas  \ref{componentlem} and \ref{lem8} and the fact that $F(w) = \infty$. 
Thus $C$ is bounded; otherwise there exist $\zeta_n \in C$ with $\zeta_n \to \infty$ and
$F( \zeta_n ) \to \zeta^* \in F(C \cup \partial C) $, contradicting the univalence of $F$ on $C$. 

Suppose that $\partial C$ has a component $\Gamma \subseteq H^+$. Then $\Gamma$ 
is a Jordan curve, and $\Gamma = \partial C$, because $C$ is simply connected.
Moreover, $\Gamma$ forms part of the boundary of  a 
multiply connected component $E$ of $W^+$. But $F$ has a pole on $\partial C$, and
$F$ is finite-valent on each such $E$, and so there are finitely many components $C$ of this type,
and none at all if $f''$ has only real zeros. 

Assume henceforth that every component of $\partial C$ meets $\R$, and take $z_0 \in \partial C$ with
the property that $ {\rm Im} \, z_0 = \max  \{ {\rm Im} \, z : \, z \in C \cup \partial C \}$.
Follow $\partial C$ in each direction, starting from $z_0$, 
until the first encounter with
$\R$. This gives 
a Jordan arc or curve $\gamma$ in $\partial C \cap
(H^+ \cup \R)$, 
such that $\gamma \cap \R = \{ a, b \}$, where $a$ and $b$ are real zeros of $F'$
with $a \leq b$. Here it is necessary to
allow for the possibility that
$a = b$, in which case $a$ is a multiple zero of $F'$ and so of $ff''$. 
Now $\lambda = \gamma \cup [a, b]$ is a Jordan curve, and since  $F'(z_0) \neq 0$ local considerations show that 
there are points in $C$ which lie in the interior domain of $\lambda$, and hence so does all of $C$.

Let $c = \sup \{ x \in \R : \, [a, x] \subseteq \partial C \}$. Then $[a, c] \subseteq  \partial C$, and $a$ and $c$ are zeros of $ff''$
(again, in principle, $a$ and $c$ might coincide, and so might $b$ and $c$). 
Lemmas \ref{lemk1} and \ref{lem11b} show that 
$f$ has
no poles in $\partial C$, each of $f$ and $f''$ has one simple zero in the set $\{ a, c \}$, and $c \leq b$.

Now $f$ has at least one pole in $(-\infty, a)$, since otherwise neither $a$ nor $c$ can be a zero of $f''$, by Lemma \ref{lem11a}. 
Let $A$ be the nearest pole of $f$ to $a$ in $(-\infty, a)$. Then $A$ lies on the boundary of a component $D$ of $W^+$. Because 
$F$ has no multiple points in $[A, a) $ by Lemma \ref{lemk1}, the interval
$[A, a] $ is a subset of $ \partial D$.
Furthermore,  $\gamma $ meets $\partial D$: if $a$ is a simple zero of $F'$ then this is clear, while
if $a$ is a multiple zero of $F'$ then $F'''(a) > 0$ by Lemma \ref{lem11cc}, in which case $\gamma $ meets $\partial D$
because $C$ lies in the interior domain of $\lambda = \gamma \cup [a, b]$. 
Since $f''$ has no non-real zeros in the closure of $C$ it follows that  $\gamma \subseteq \partial D$. 
A similar argument shows that there exists a pole $B$ of $f$ with $B > b$, such that the interval $[b, B]$ lies in the boundary of 
a component $D'$ of $W^+$, and so does $\gamma$, from which it follows that $D = D' =D_w$.


In the case where $f''$ has only real zeros, $F$ must be univalent on $D$, and the branch $g$ of the inverse function $F^{-1}$ which maps $H^+$
to $D$ has at least two attracting fixpoints on the boundary of $ H^+$, at $A$ and $B$, contradicting the Denjoy-Wolff theorem 
\cite[Chapter 2]{Stei2}. Indeed, the iterates $g^n$ form a normal family on $H^+$, since $g(H^+) = D \subseteq H^+$, but 
$g$ extends to be analytic on a neighbourhood $U_A$ of $A$, such that $g(U_A) \subseteq U_A$
and the $g^n$ converge to $A$ on $U_A$, and in the same way they 
converge to $B$ on a neighbourhood of $B$. 

In the general case where $f''$ has finitely many non-real zeros, suppose that there exist infinitely many zeros $w \in H^+$ of $f'$.
This gives infinitely many distinct components $C_w$ of $W^-$
as above, each with a corresponding component $D_w$ of $W^+$. The $D_w$ need not be distinct,
but Lemma \ref{lemvalency} implies that $L$ has finitely many poles on the boundary of any component of $W^+$, and therefore so has
$f$. Hence there must exist at least one $D_w$ which is mapped univalently onto $H^+$ by $F$, and the Denjoy-Wolff theorem 
supplies a contradiction as before. 
\hfill$\Box$
\vspace{.1in}

To complete the proof of Theorem \ref{kreinthm}, it now follows 
from Lemma \ref{lem11} and the fact that all but finitely many zeros of $f$ and $f''$ are real that $f$ satisfies the hypotheses of
\cite[Theorem 6.4]{ankara} (see also \cite[Theorem 1.5]{Laams09}), subject to the assumption made earlier that $A = B =0$ in (\ref{krein2}).
Then 
$$
 f(z) = \frac{ R(z)e^{icz} -1}{A_1 R(z) e^{icz} - \overline{A_1}} \, , 
$$
with  $c \in (0, \infty) $, $ A_1 \in \C \setminus \R$, 
and $R$ a rational function satisfying 
$|R(x)| = 1$ for all $x \in \R$, by \cite[Theorem 6.4]{ankara}. 
Since all residues of $f$ have to be negative, it follows easily that $A_1 \in H^+$, 
and the fact that $f(H^+) \subseteq H^+$ shows that all zeros of $R$ lie in $H^+$, and all poles in $H^-$.

Finally, suppose that all zeros of $f''$ are real. Then the Schwarzian derivative $S_f$ is entire, because $f'$ has no zeros and 
all poles of $f$ are simple  \cite{Hil1,Hil2}. 
Since $f$ is transcendental of order at most~$1$, it must be the case that $S_f$ is a non-zero constant, so that there exist $a  \in \C$ and
a M\"obius transformation $T$ such that $f(z) = T( e^{i 2 a z} )$. Because $f$ is real with only real zeros and poles, $a$ must be real,
and $f(z) = C \tan (az+b) + E$, with  $b$, $C$ and $E$ also real. 
\hfill$\Box$
\vspace{.1in}

\section{A special case of Theorem \ref{f''=0} }

The following special case illustrates Theorem \ref{f''=0} and plays a key role in its proof. 
\begin{lem}
 \label{lemresidue}
Let $a, b, D, E \in \C$ with $a \neq 0$ and $D \neq E$, and let $2 \leq n \in \Z$. Let 
\begin{equation}
 \label{l0}
F(z) = \left( \frac{De^{az+b} - E}{e^{az+b}-1} \right)^n.
\end{equation}
(i) There exists a meromorphic function $G$ in the plane with $G' = F$ if and only if $D = \lambda E$ 
where $\lambda^n = 1$, $\lambda \neq 1$. \\
(ii) There does not exist a meromorphic function $H$ in the plane with $H'' = F$. 
\end{lem}
\textit{Proof.} It may be assumed that $a=1$ and $b=0$.
By periodicity, there exists a meromorphic function $G$ with $G' = F$ if and only if ${\rm Res} \, (F, 0) = 0$. 
The function $w = e^z - 1 $ 
is univalent on a neighbourhood of the origin and has local inverse 
\begin{equation}
\label{l2} 
z = \phi(w) = \log (1+w) = w - \frac{w^2}2 + \frac{w^3}3 - \ldots .
\end{equation}
Let $\varepsilon$ be small and positive and let $\gamma$ describe the circle $|z| = \varepsilon$ once 
counter-clockwise. Let $\Gamma$ be the image of $\gamma$ under $w = e^z - 1$. Then ${\rm Res} \, (F, 0) = 0$
if and only if 
\begin{equation}
 \label{l3}
0 = \int_\gamma F(z) \, dz = \int_\Gamma \psi(w) \, dw, \quad 
\psi (w) = \left( D + \frac{D-E}w \right)^n \phi'(w) .
\end{equation}
Now (\ref{l2}) and (\ref{l3}) give, as $w \to 0$,  
$$
\psi(w) = \left( D^n + n D^{n-1} \left( \frac{D-E}w \right) + \ldots + \left( \frac{D-E}w \right)^n \right) 
\left(1 - w + \ldots + (-1)^{n-1} w^{n-1} + \ldots \right) ,
$$ 
and so (i) follows from the fact that
\begin{eqnarray*}
 {\rm Res} \, ( \psi, 0) &=& 
nD^{n-1} (D-E) - \frac{n!}{2!(n-2)!} D^{n-2} (D-E)^2 + \ldots + (-1)^{n-1} (D-E)^n \\
&=& - \left( nD^{n-1} (E-D) + \frac{n!}{2!(n-2)!} D^{n-2} (E-D)^2 + \ldots +  (E-D)^n \right) \\
&=& - \left( (D + E-D)^n - D^n \right) = D^n - E^n. 
\end{eqnarray*}

To establish (ii), suppose that there does exist a meromorphic function $H$ in the plane with $H'' = F$. Then 
$D = \lambda E$, with $\lambda^n = 1$  by (i), and it may be assumed that $E=1$ and $D = \lambda \neq 1$. This time write
\begin{equation}
 \label{l4}
w = q(z) =  \frac{e^z-1}{\lambda e^z - 1} , \quad z = q^{-1}(w) = \sigma(w) = \log \left( \frac{1-w}{1-\lambda w} \right) ,
\end{equation}
each of these being univalent near the origin. 
This forces, with $\gamma$ as before and $\Lambda $ the image of $\gamma$ under $w = q(z)$,
\begin{equation}
 \label{l5}
0 =  \int_\gamma z F(z) \, dz =  \int_\gamma \frac{z}{w^n}  \, dz = \int_\Lambda \frac{ \tau (w) }{w^n} \, dw ,
\quad \tau (w) = \sigma (w) \sigma'(w).
\end{equation}
Now, as $w \to 0$, expanding (\ref{l4}) yields
\begin{eqnarray*}
\label{l6}
 \tau (w) &=& 
\left( w(\lambda-1) + \ldots + \frac{w^{n-1}}{n-1} \left( \lambda^{n-1}-1\right) + \ldots \right)  
\left( \lambda - 1 + \ldots + w^{n-2} \left( \lambda^{n-1} -1 \right) + \ldots \right) \nonumber \\
&=& a_1 w + \ldots + a_{n-1} w^{n-1} + \ldots . . 
\end{eqnarray*}
Here the coefficient $a_{n-1}$ of $w^{n-1}$ must vanish by (\ref{l5}), which delivers
\begin{equation}
 \label{l7}
0 = \frac1{n-1} \left( \lambda^{n-1}-1\right) ( \lambda - 1) + \ldots + ( \lambda - 1)\left( \lambda^{n-1} -1 \right) 
= \sum_{j=1}^{n-1} \frac1{n-j} \left( \lambda^{n-j}-1\right)\left( \lambda^{j}-1\right) .
\end{equation}
But $\lambda^n = 1$ and so $\lambda = \exp( 2 \pi i k/n )$ for some $k \in \{ 1, \ldots , n-1 \}$. It follows that,
for $1 \leq j \leq n-1$, 
$$
\mu_j = \left( \lambda^{n-j}-1\right)\left( \lambda^{j}-1\right) = 2 - \left( \lambda^j + \lambda^{-j} \right) =
2 - 2 \cos ( 2 \pi jk/n ) \geq 0.
$$
Since $\mu_1 > 0$, the sum in (\ref{l7}) is real and positive, and this contradiction completes the proof.
\hfill$\Box$
\vspace{.1in} 
\section{Proof of Theorem \ref{f''=0}}
Let $f$ be as in the hypotheses, let $R$ be a large positive real number, and define $g$ formally by
\begin{equation}
 \label{1}
f' = g^n .
\end{equation}
Then $g$ admits unrestricted analytic continuation in $R < |z| < \infty$, these continuations having only simple poles and
no critical points. Since $g'/g$ is single-valued in the plane, so is the function $A$ defined by
\begin{equation}
 \label{2}
2A = S_g = \frac{g'''}{g'} - \frac32 \left( \frac{g''}{g'} \right)^2 ,
\end{equation}
where $S_g$ denotes the Schwarzian derivative \cite{Hil1,Hil2}. 
Moreover, $A$ has finitely many poles, and none in $R < |z| < \infty$,
because the continuations of $g$  are free of multiple
points there. 
\begin{lem}
 \label{lemnotrsp}
The function $A$ is rational but does not satisfy $A(z) = O( |z|^{-2})$ as $z \to \infty$.
\end{lem}
\textit{Proof.} 
The first assertion follows from the lemma of the 
logarithmic derivative and the fact that $f$ has finite lower order. 
Now suppose that $A(z) = O( |z|^{-2})$ as $z \to \infty$. Take $z_0 \in \C$ with $|z_0| > R$ 
such that $z_0$ is neither a pole nor a zero of $f'$, and define the functions $W$ and $V$ in a simply connected open neighbourhood $U$
of $z_0$ by 
\begin{equation}
 \label{3}
W^2 = \frac1{g'} = \frac{n g^{n-1}}{f''}, \quad 
V = W^{2n} = \frac{n^n (f')^{n-1}}{(f'')^n} .
\end{equation}
It follows from (\ref{3}), hypothesis (ii) and the fact that $R$ is large that $V$ extends to be analytic in
$R < |z| < \infty$, with a zero of multiplicity $2n$ at each pole of $f$, and no other zeros. In particular, $V$ has an essential singularity
at infinity. By a result of Valiron \cite[p.15]{Valiron}, the function $V$ may be written in the form 
\begin{equation}
 \label{valiron}
V(z) = z^q Y(z) (1+o(1))  \quad \hbox{as $z \to \infty$,}
\end{equation}
in which $q$ is an integer and $Y$ is a transcendental entire function.

A standard calculation starting from (\ref{2}) and (\ref{3}) shows that $W$ is a solution on $U$ of 
\begin{equation}
 \label{4}
w'' + A(z) w = 0.
\end{equation}
On the other hand, (\ref{3}) and (\ref{4}) now yield, again on $U$,  
\begin{equation}
 \label{5}
W = V^{1/2n} , \quad - A = \frac{W''}{W} , \quad 
-A =  \frac1{2n} \left( \frac1{2n} - 1 \right) \left( \frac{V'}{V} \right)^2 + \frac1{2n} \frac{V''}{V} .
\end{equation}
The last equation of (\ref{5}) then holds by analytic continuation throughout $R < |z| < \infty$. 

Now let $\nu(r)$ denote the central index of the transcendental entire function $Y$. By (\ref{valiron}) and 
the Wiman-Valiron theory \cite{Hay5}, if $r$ is large and 
lies outside a set of finite logarithmic measure, and if $|z_1| = r$ and $|Y(z_1)| = M(r, Y)$, then $\nu(r)$ is large and 
$$
\frac{V'(z_1)^2}{V(z_1)^2} \sim \frac{V''(z_1)}{V(z_1)} \sim \frac{\nu(r)^2}{z_1^2}  \quad \hbox{and} \quad 
 \frac1{4n^2} \frac{\nu(r)^2}{z_1^2} \sim - A(z_1) = O( r^{-2} ),
$$
which is a contradiction.
\hfill$\Box$
\vspace{.1in}
 
Lemma \ref{lemnotrsp} makes it possible to write, as $z \to \infty$, 
\begin{equation}
 \label{6}
A(z) \sim c z^m , \quad c \in \C \setminus \{ 0 \}, \quad m \in \Z , \quad m \geq -1, 
\end{equation}
and so Hille's asymptotic method \cite{Hil1,Hil2} may now be applied to (\ref{4}). 
The $m+2$ critical rays $\arg z = \theta_0 $ for the equation (\ref{4}) are determined by the formula
\begin{equation}
 \label{7}
\arg c + (m+2) \theta_0 = 0 \quad \hbox{ ( mod $2 \pi$).} 
\end{equation}
Let $\varepsilon$ and $1/R_1$ be small and positive: then (\ref{4}) has linearly independent 
solutions $u_1, u_2$ satisfying
\begin{equation}
 \label{9}
u_1(z) \sim A(z)^{-1/4} e^{-iZ} ,  \quad u_2(z) \sim A(z)^{-1/4} e^{iZ} , \quad Z = \int_{2R_1}^z A(t)^{1/2} \, dt \sim \frac{2c^{1/2}}{m+2} z^{(m+2)/2} ,
\end{equation}
as $z \to \infty$
in the sectorial region
$$
S (R_1,  \varepsilon ) = \left\{ z \in \C : |z| > R_1 , \, | \arg z - \theta_0 | < \frac{2 \pi}{m+2} - \varepsilon \right\} .
$$
If $m=-1$ then there is only one critical ray given by (\ref{7}), and $S(R_1, \varepsilon ) $ 
should be understood as
lying on the Riemann surface of $\log z$. It follows from (\ref{1}), (\ref{2}) and (\ref{4}) that there exist 
complex numbers $A_j$ and $B_j$ such that $f'$ satisfies, on $S (R_1,\varepsilon )$, 
\begin{equation}
 \label{10}
f' = g^n , \quad g = \frac{ A_1 u_1 - A_2 u_2}{B_1 u_1 - B_2 u_2}  ,
\end{equation}
and $A_1B_2 - A_2 B_1 \neq 0$, since $f'$ is non-constant. 

It may be assumed that $\theta_0$ is chosen so that $f$ has infinitely many poles in 
the narrower sectorial region $S(R_1,  4 \varepsilon )$, 
which forces $B_1 B_2 \neq 0$ in (\ref{10}) and makes it possible to write
\begin{equation}
\label{11}
 f' = \left(  \frac{D e^{2 \pi i L}  - E}{e^{2 \pi i L} - 1 } \right)^n , \quad D, E \in \C, \quad D \neq E ,
\end{equation}
where 
\begin{equation}
 \label{12}
L(z) =  \frac{  1}{2 \pi i} \log \left( \frac{B_2u_2(z)}{B_1u_1(z)} \right)  \sim  \frac{ Z}{\pi}  \sim \frac{2c^{1/2}}{\pi (m+2) } z^{(m+2)/2}
\end{equation}
as $z \to \infty$ in $S(R_1,  2 \varepsilon )$. In view of (\ref{7}) it may be assumed that 
the branch of the square root in (\ref{9}) is chosen so
as to make ${\rm Re} \, L(z)$ positive as $z \to \infty$ on the critical ray, and 
the poles $\zeta_j$  of $f$ in $S(R_1, 4 \varepsilon )$ must have $\arg \zeta_j \to \theta_0$ as $\zeta_j \to \infty$. 

The asymptotics (\ref{12}) show that $w = L(z)$ maps a subdomain $S^*$ of $S(R_1,  3 \varepsilon )$ univalently onto a 
a sectorial region $\Omega = \{ w \in \C : |w| > R_2, \, | \arg w | < \pi - \delta \}$, 
where $R_2$ is large, and $\delta$ may be made arbitrarily small by choosing $\varepsilon$ small enough. 
In particular, $\Omega$ contains a half-plane 
$H$ given by ${\rm Re} \, w > q_0 > 0$. 
Let $z = \phi (w)$ be the inverse mapping from $\Omega$ to $S^*$, choose a large positive integer $q$ and 
let the contour $\gamma$ in $H$ describe once counter-clockwise the circle of centre $q$ and radius $1/4$. 
Then $f$ has no poles on $\phi(\gamma)$ and (\ref{11}) gives
\begin{equation}
 \label{13}
0 = \int_{\phi(\gamma)} f'(z) \, dz = \int_\gamma \psi(w) \, dw, \quad \psi (w) =
\left( \frac{D e^{2 \pi i w} - E}{e^{2 \pi i w} - 1 } \right)^n \phi'(w) .
\end{equation}
As $w \to q$ periodicity yields
$$
Q(w)
= \left(  \frac{D e^{2 \pi i w} - E}{e^{2 \pi i w} - 1 } \right)^n = \left(  \frac{D e^{2 \pi i (w-q)} - E}{e^{2 \pi i (w-q) } - 1 } \right)^n =
\frac{D_n}{(w-q)^n} + \ldots + \frac{D_1}{w-q} + O(1) ,
$$
in which the $D_j$ depend on $n$, $D$ and $E$ but not on $q$.
Moreover, Lemma \ref{lemresidue} implies that the function $Q(w)$ 
is not the second derivative of a meromorphic function in the plane and so, by periodicity again, at least one of $D_1$ and $D_2$ is non-zero. 
Now  (\ref{13}) delivers
\begin{equation*}
 \label{14}
0 = {\rm Res} \, ( \psi, q) = \sigma (q), \quad
\sigma (w) = D_1 \phi'(w) + D_2 \phi''(w) + \ldots + D_n \frac{\phi^{(n)}(w)}{(n-1)!} , \quad |D_1| + |D_2| > 0 . 
\end{equation*}
Since $m+2 \geq 1$ in (\ref{12}), the function $\sigma (w)$ has at most polynomial growth in the half-plane ${\rm Re} \, w > q_0 + 1$.
Now the fact that  $\sigma (q) = 0$ for all sufficiently large positive integers $q$ forces
$\sigma $ to vanish identically (using, for example, \cite[Lemma 5]{La9}). 
This implies that $\phi$ satisfies, in the domain $\Omega$,
a linear differential equation with constant
coefficients, and so $\phi$ is an entire function of exponential type. Because $\phi$ has polynomial growth in $\Omega$, by (\ref{12}),
while $\delta $ is small, applying
the Phragm\'en-Lindel\"of principle  shows that $\phi$ is a polynomial.
But then the condition $|D_1| + |D_2| > 0$ and the vanishing of $\sigma$ together ensure that  
$\phi$ is a polynomial of degree $1$, and so is its inverse function~$L$. 
Thus (\ref{11}) implies that Lemma \ref{lemresidue} may be applied to $f'$, which completes the proof.  
\hfill$\Box$
\vspace{.1in}

\section{Proof of Theorem \ref{f''=0real}}\label{pfthmreal}

Let $f$ be a real transcendental meromorphic function in the plane 
satisfying hypotheses (a), (b) and (c) of Theorem \ref{f''=0real}. 
It is not assumed at this stage that hypothesis (d) holds.   
The function
\begin{equation}
 \label{w1}
h = \frac{f'}{f''} 
\end{equation}
has finitely many poles and non-real zeros. If $h$ is a rational function
then $f' = R_0 e^{P_0}$ with $R_0$ a real rational function and $P_0$ 
a real polynomial. Because $f$ has finitely many non-real zeros, this forces (\ref{ReP}). 
Assume for the remainder of the proof that $h$ is transcendental.

\begin{lem}
 \label{lemw1}
The function $L = f'/f$ is transcendental, and its Tsuji characteristic satisfies $T_0(r, L) = O( \log r )$ 
as $r \to \infty$. 
\end{lem}
\textit{Proof.} $L$ must be transcendental, because $1/h = L + L'/L$. The second assertion holds by Lemma \ref{lem1}
and the fact that all but finitely many zeros and poles
of $f$ and $f''$ are real. 
\hfill$\Box$
\vspace{.1in}

\begin{lem}
 \label{lemw2}
The Nevanlinna characteristic of $h$ satisfies $T(r, h) = O( r \log r )$ as $r \to \infty$, while
\begin{equation}
 \label{w2}
\overline{N} (r, f) + \overline{N} (r, 1/f)  + \overline{N} (r, 1/f') = O( r \log r ) \quad \hbox{as $r \to \infty$.}
 \end{equation}
Furthermore, $T(r, L) = O( r \log r )$ as $r \to \infty$.
\end{lem}
\textit{Proof.} Lemma \ref{lemw1} and standard properties of the Tsuji characteristic give $T_0(r, h) = O( \log r )$ 
as $r \to \infty$, so that $T(r, h) = O( r \log r )$ as $r \to \infty$ by Lemma \ref{lem2}.
It then follows that
$$
\overline{n} (r, f) +  \overline{n} (r, 1/f') \leq \overline{n}(r, 1/h) = O( r \log r ) \quad \hbox{as $r \to \infty$,}
$$
using (\ref{w1}). The corresponding result for $\overline{n} (r, 1/f)$ now follows from Rolle's theorem. This gives (\ref{w2}),
which implies the estimate
for $T(r, L)$, using Lemmas \ref{lem2} and \ref{lemw1}.
\hfill$\Box$
\vspace{.1in}


\begin{lem}
 \label{lemw4}
The function $f$ admits a representation
\begin{equation}
 \label{w3}
f = \frac{G}{H} , \quad \frac{G'}{G} = \phi \psi ,
\end{equation}
in which:\\
(i) $G$ and $H$ are real entire functions, and $H$ has order at most $1$;\\
(ii) $\phi$ and $\psi$ are real meromorphic functions, and $\phi$ has finitely many poles and order at most $1$;\\
(iii) either $\psi \equiv 1$ or $\psi$ maps the upper half-plane $H^+$ into itself. 
\end{lem}
\textit{Proof.} Here $H$ is the canonical product formed using the poles of $f$, all but finitely many of which are real, the rest
occurring in conjugate pairs because $f$ is real. Since the poles of $f$ have bounded multiplicities, it 
follows from (\ref{w2}) that $H$ has order at most $1$. Now $G$ is a real entire function with finitely many non-real zeros,
and the formula $G'/G = \phi \psi$ is just the standard Levin-Ostrovskii factorisation \cite{BEL,lajda}, 
in which $\psi$ is formed as in the proof of Lemma \ref{lem4}, using real zeros $a_k$ of $G$ and $b_k$ of $G'$. Finally,
$\phi$ has order at most $1$ because (\ref{C1}) holds with $f$ replaced by $\psi$ so that, as $r \to \infty$, 
$$
m(r, \phi) \leq m(r, G'/G) + m(r, 1/\psi)   \leq m(r, G'/G) + O( \log r) \leq m(r, L) + O( \log r) .
$$
\hfill$\Box$
\vspace{.1in}


\begin{lem}
 \label{lemw5}
The function $\phi$ in (\ref{w3}) is rational, and $G$ and $f$ have finite order. 
\end{lem}
\textit{Proof.} Assume that $\phi$ is transcendental. 
Fix  a small positive real number $\varepsilon$ and a large positive integer $N$, and set 
\begin{equation}
 \label{w4}
W_1(z) = \frac{h(z)}{z^N} = \frac{f'(z)}{z^N f''(z)} , \quad W_2(z) = \frac{\phi(z)}{z^N} .
\end{equation}
Each $W_j$ has finite order and finitely many poles, and so  Lemma \ref{lemcarleman}
gives an unbounded set $E_1 \subseteq [1, \infty )$ such that for $r \in E_1$ 
and $j=1, 2$ there exists $\theta_j \in \R$ with 
\begin{equation}
 \label{w5}
|W_j(re^{i\theta} ) | \geq 1 \quad \hbox{for} \quad | \theta - \theta_j | \leq 8 \varepsilon .
\end{equation}
For $r \in E_1$, integration gives $c_r \in \C \setminus \{ 0 \}$ 
and $d_r \in \C$ such that 
$$
f'(re^{i\theta}) = c_r  \left( 1 + O  \left( r^{1-N} \right) \right), \quad 
f(re^{i\theta}) = c_r  \left( re^{i\theta}  +  O \left( r^{2-N} \right) \right) + d_r  
$$
for $| \theta - \theta_1 | \leq 8 \varepsilon $. This gives in turn,
for $\theta$ in an interval of length $4 \varepsilon$, 
\begin{equation}
 \label{w6}
P( re^{i\theta} )
=re^{i\theta} \frac{f'( re^{i\theta} ) }{f ( re^{i\theta} ) } = \frac{ re^{i\theta} (1+o(1))}{ re^{i\theta} + d_r/c_r + o(1)} =
O(1) .
\end{equation}
Because $f$ is real it may be assumed that (\ref{w6}) holds for at least one $\theta$ in the 
interval $[ \varepsilon, \pi - \varepsilon ]$, and so Lemma \ref{norfamlem} yields $P( re^{i\theta} ) = O(1)$
for $r \in E_1 $ and all $\theta \in [ \varepsilon, \pi - \varepsilon ]$. Since $H$ has  order at most $1$
and finitely many non-real zeros,
(\ref{C1}), with $f$ replaced by $\psi$, and (\ref{w3}) yield 
$$
\frac{G'( re^{i\theta} ) }{G ( re^{i\theta} ) } = \frac{f'( re^{i\theta} ) }{f ( re^{i\theta} ) } +
\frac{H'( re^{i\theta} ) }{H ( re^{i\theta} ) } = O(r) \quad \hbox{and} \quad 
\phi ( re^{i\theta} ) =  \frac{G'( re^{i\theta} ) }{G ( re^{i\theta} ) \psi ( re^{i\theta} )} = O(r^2) 
$$
for $r \in E_1 $ and $|\theta | \in [ \varepsilon, \pi - \varepsilon ]$. By (\ref{w4}) this contradicts (\ref{w5}) for $j=2$.

Thus $\phi$ is rational, and the assertion that $G$ has finite order, which in turn implies that so has $f$, follows from
a standard argument \cite[Lemma 5.1]{BEL}. 
\hfill$\Box$
\vspace{.1in}

\begin{lem}
 \label{lemw5a}
The function $f'$ has finitely many asymptotic values, all transcendental singularities of the inverse function of $f'$ are
logarithmic, and $f''/f'$ has lower order at least $1/2$. 
\end{lem}
\textit{Proof.} Since $f''/f'$ has finitely many zeros, $f'$ has finitely many critical values. Thus,
because $f'$ has finite order, all transcendental singularities
of the inverse function are direct, by the main result of \cite{BE}, and they are finite in number by the 
Denjoy-Carleman-Ahlfors theorem \cite{Hay7}. Hence all such singularities are in fact logarithmic. 

The last assertion is proved as in \cite[Lemma 11]{Lagehring}. 
Since $f''/f'$ has finitely many zeros, the same result of Lewis, Rossi and Weitsman \cite{LRW} as used in Lemma \ref{lem5}
gives a path $\gamma$ tending to infinity on which $f'$ tends to $\beta \in \C \setminus \{ 0 \}$. 
If $f''/f'$ has lower order less than $1/2$ then the $\cos \pi \rho$ theorem \cite{Hay7} implies that $f''/f'$ is small, and 
$f'$ is close to $\beta$, on the union of a sequence of circles $|z| = r_n \to \infty$. This contradicts the fact 
that the singularity over $\beta$ is logarithmic. 
\hfill$\Box$
\vspace{.1in}

\begin{lem}
 \label{lemw6a}
Let $\delta_1 > 0$ and let $\rho < \infty$ be the order of growth of $f$. Then 
$\left| f''(z)/f'(z) \right| \leq |z|^{\rho} $
as $z \to \infty$ with $\delta_1 \leq | \arg z | \leq \pi - \delta_1 $. 
\end{lem}
\textit{Proof.} This follows from standard estimates based on the differentiated Poisson-Jensen formula~\cite{Hay2}
and the fact that $f'$ has order  $\rho$ and finitely many non-real zeros and poles.
\hfill$\Box$
\vspace{.1in}

\begin{lem}
 \label{lemw7a}
There exists $\alpha \in \C \setminus \{ 0 \}$ with the following property. 
If $\varepsilon > 0$ then, 
as $z \to \infty$ with $\varepsilon \leq  \arg z  \leq \pi - \varepsilon $, 
\begin{equation}
 \label{w8a}
\left| \frac{f''(z)}{f'(z)} \right| \leq \exp \left( -|z|^{1/4} \right) 
\end{equation}
and $f'(z) = \alpha + o(1)$. 
\end{lem}
\textit{Proof.} To prove (\ref{w8a}) apply Lemma \ref{lemhfn} with $Q= f''/f'$ and $\eta = 1/16$, in conjunction with Lemmas
\ref{lemw5a} and \ref{lemw6a}.  
Integration then gives $f'(z) = \alpha + o(1)$ in the same sector, where $\alpha \in \C \setminus \{ 0 \}$,
and it is clear that $\alpha$ is independent
of $\varepsilon$. 
\hfill$\Box$
\vspace{.1in}

\begin{lem}
 \label{lemw8a}
The inverse function of $f'$ has exactly one of the following:\\
(I) a logarithmic singularity over each of $\alpha$ and $\bar \alpha$, where $\alpha \in \C \setminus \R$, and
no other transcendental singularities;\\
(II) one or two logarithmic singularities over $\alpha \in \R \setminus \{ 0 \}$, and no other transcendental singularities. 
\end{lem}
\textit{Proof.}
Lemma \ref{lemw7a} gives $f'(z) = \bar \alpha + o(1) $ 
as $z \to \infty$ with $\varepsilon \leq  - \arg z  \leq \pi - \varepsilon $,
where $\varepsilon $ may be chosen arbitrarily small. 
The result now follows from Lemmas \ref{lemcarleman} and  \ref{nodirect}. 
\hfill$\Box$
\vspace{.1in}

Following \cite{Lagehring}, let $J$ be a  polygonal Jordan curve in $\C \setminus \{ 0 \}$,
symmetric with respect to the real axis, such that every
finite non-zero critical or asymptotic value of $f'$ 
lies on $J$ but is not a vertex of $J$.
Here $J$ can be formed so that 
its complement in $\C \cup \{ \infty \} $
consists of two simply connected domains $B_1$ and $B_2$, with
$0 \in B_1$ and $\infty \in B_2$. Fix conformal
mappings 
\begin{equation}
h_m : B_m \to  \{ w \in \C : |w| < 1 \}, 
\quad m = 1, 2, \quad h_1(0)=0, \quad h_2( \infty ) = 0.
\label{a1}
\end{equation}
The mapping $h_1$ may then be extended to be quasiconformal on the plane \cite[Ch.5]{Pom2}, fixing infinity, and there
exist a meromorphic function $G_1$  and a quasiconformal mapping $\psi_1$ such that 
\begin{equation}
 \label{G1def}
h_1 \circ f' = G_1 \circ \psi_1 \quad \hbox{on $\C$. }
\end{equation}
The following is \cite[Lemma 4]{Lagehring}, translated to the present setting in the light of Lemma \ref{lemw8a}.

\begin{lem}
 \label{complem}
For $j=1,2$, all components of $(f')^{-1}(B_j)$ are simply connected  and all but finitely
many  are unbounded. If $C_0$ is a component  of $(f')^{-1}(B_1)$ then 
$C_0$ contains one zero of $f'$,  of multiplicity $m_1 \in \N$, and $C_0$ is mapped $m_1$ to $1$ onto $ B_1$ by $f'$.
Furthermore, if a zero $z_1$ of $f''$ lies in a component $C_1$ of $(f')^{-1}(B_1)$ then 
$z_1$ is the only zero of $f''$ in~$C_1$. 
Similarly, each component of $(f')^{-1}(B_2)$ contains exactly one pole of $f$, disregarding multiplicities.
\end{lem}
\hfill$\Box$
\vspace{.1in}

The next step is to combine  \cite[Lemma 5]{Lagehring}  with Lemma \ref{lemw8a}. 

\begin{lem}\label{lem3a}
Arbitrarily small positive real numbers
$\varepsilon_1$ and $\varepsilon_2$ may be chosen with the following properties.
There exist one or two
unbounded simply connected domains $U_n$, each  a component of the set
$\{ z \in \C : | f'(z) - b_n | < \varepsilon_1 \}$,
such that $U_n$ 
contains a
path tending to infinity on which $f'(z)$ tends to  $b_n $. Here each
$b_n$ is $\alpha$ or $\bar \alpha$, and 
$f'(z) \not = b_n$ on $U_n$, while 
$|f(z) - b_n z| < \varepsilon_2 |z|$ for all
$z$ in $U_n$  with $|z|$ large enough.
If $\Gamma$ is a path tending to infinity on which
$f'$ tends to an asymptotic value $\beta$, then there exists $n$ such that
$\beta = b_n$ and $\Gamma \setminus U_n$ is bounded. 
\end{lem}
\hfill$\Box$
\vspace{.1in}

\begin{lem}
 \label{lemboundary}
The function $f'$ has infinitely many zeros $x_j$, all but finitely many of which satisfy the following. 
First, $x_j$ is real and lies in a  component $C_j$ of
$(f')^{-1}(B_1) $ which  is unbounded, simply connected and symmetric with respect to the real axis, 
and there are no zeros of $f''$ on the boundary $\partial C_j$. Furthermore,  $\partial C_j $ is  $ \Gamma_j^- \cup \Gamma_j^+$,
where each $\Gamma_j^\pm$ is a simple curve tending to infinity
in both directions, symmetric with respect to $\R$,
and meeting the real axis exactly once. Analogous considerations apply to poles of $f'$. 
\end{lem}
\textit{Proof.}
There exist infinitely many zeros $x_j$ of $f'$ by Lemma \ref{lemw8a}. For $|x_j|$ large let  
\begin{equation}
 \label{y1}
 -  \infty <  y_j^- = \inf \{ x \in \R: [x, x_j] \subseteq C_j \} < 
y_j^+ = \sup \{ x \in \R: [x_j, x] \subseteq C_j \} < \infty . 
\end{equation}
Each $y_j^\pm$ lies in a component $\Gamma_j^\pm$ of $\partial C_j$ which is symmetric with respect to $\R$, and
$\psi_1(\Gamma_j^\pm ) $ is a level curve of the function $G_1$ in (\ref{G1def}). Thus
$\Gamma_j^\pm \cap \R = \{ y_j^\pm \}$, because $C_j$ is simply connected. 
Finally, observe that any component of 
$\partial C_j$ other than the $\Gamma_j^\pm $
would have to lie in $\C \setminus \R$ and form part of the boundary of a component of $(f')^{-1}(B_2)$, 
that component having to contain a non-real pole of $f$. 
\hfill$\Box$
\vspace{.1in}

\begin{lem}
 \label{lemy1}
The zeros of $f'$ have bounded multiplicities, and case (II) holds in Lemma \ref{lemw8a}. 
\end{lem}
\textit{Proof.}
Each $\Gamma_j^\pm$ in Lemma \ref{lemboundary}
forms part of the boundary of a a component of $(f')^{-1}(B_2)$, and  the poles of $f$ have
bounded multiplicities. Hence the variation of $\arg f'$ on $\Gamma_j^\pm$  has an upper bound which is independent of $j$,
thus proving the first assertion.

Suppose now that case (I) holds in Lemma \ref{lemw8a}. If $z_0$ is large and is a zero of 
$f''$ then $z_0 $ and $f(z_0)$ are real, so that 
$$| f(z_0) - \alpha z_0 | = | f(z_0) - \bar \alpha z_0 | \geq  |z_0 \, {\rm Im} \, \alpha \, |.$$ 
Theorem \ref{theorem1} may now be applied, to conclude that
$f'' = R_2 e^{P_2}$ with $R_2$ a real rational function and $P_2$ a real polynomial. 
Thus $f$ has finitely many poles, which contradicts  Lemma \ref{lemw8a}.
\hfill$\Box$
\vspace{.1in}

It may be assumed henceforth that case (II)  holds in Lemma \ref{lemw8a}, with $\alpha = 1$. 

\begin{lem}
 \label{lemy2}
Fix positive real numbers $M_1$ and $M_2$ with 
$M_1$ large and $M_1 < M_2$. 
Let $v_j \in \R$ with $|v_j|$ large be a pole of $f$ of multiplicity $m_j$, and let $D_j$ be the component of 
$(f')^{-1} (B_2)$ in which $x_j$ lies. Then $|f(z) - z| \leq 2 \varepsilon_2 |z|$ for all $z \in D_j$ with $M_1 < |f'(z)| < M_2$, 
where $\varepsilon_2$ is as in Lemma \ref{lem3a}. Moreover,
$f$ has at least $m_j$ real simple zeros in $D_j$, 
and 
$m_j$ is $1$ or $2$. 
\end{lem}
\textit{Proof.}
The component $D_j$ is simply connected and, as shown in
Lemma \ref{lemboundary}, its boundary consists of two disjoint simple curves $\Lambda_j^\pm$.
The function 
$v = (h_2 \circ f')^{1/m_j}$ maps $D_j$ conformally onto the unit disc, 
and as $z$ tends to infinity in either direction along either of the $\Lambda_j^\pm$, the image $f'(z)$ tends to 
the unique asymptotic value $1$ of $f'$, 
since $f'$ is finite-valent on $D_j$. This implies that $D_j$ meets one of the components $U_n$ of Lemma \ref{lem3a}. 
It follows that there exist $\mu_j$ with $ \mu_j^{m_j}  = h_2(1)$ and a  positive $\varepsilon_3$ such that
if $z \in D_j$ and  
$|v(z)-\mu_j| \leq \varepsilon_3$ then $z \in U_n$. Here $\varepsilon_3$ may be chosen arbitrarily small and 
independent of $j$, since the $m_j$ are bounded by hypothesis. 

Let $u$ be the inverse function of $v$, mapping the unit disc onto $D_j$. Then 
$u'(0) = o( |v_j| )$, by 
Koebe's $1/4$ theorem and 
Lemma \ref{lemw7a}. Koebe's distortion theorem then yields 
$u'(w) = o( |v_j| )$ for $|w| \leq 1 - \varepsilon_3$. Now let $z_1 \in D_j$ be such that
$w_1 = v(z_1)$ satisfies $\varepsilon_3 \leq |w_1| \leq 1- \varepsilon_3$. Then 
$w_1$ can be joined to a point $w_2$ with $|w_2| < 1$, $|w_2 - \mu_j| \leq \varepsilon_3$ by a path $\Sigma$ in
$\varepsilon_3 \leq |w| \leq 1- \varepsilon_3$ so that $\sigma = v(\Sigma)$ 
is a path in $D_j$, of length $o(|v_j|)$, joining $z_1$ to  $z_2 = u(w_2) \in U_n$.
But then $|f(z_2) - z_2 | \leq \varepsilon_2 |z_2|$ by Lemma \ref{lem3a}. Since $f'$ is bounded on $\sigma$, 
integration of $f'$ gives $|f(z_1) - z_1 | \leq 2 \varepsilon_2 |z_1|$, proving the first assertion.

Next, let $\tau$ be the image under $u$ of the circle $|w| = \varepsilon_3$. Then $\tau$ is a Jordan curve in $D_j$  enclosing
$v_j$, and symmetric with respect to the real axis. Furthermore, $|f(z) - z| < |z|$ on $\tau$; thus
Rouch\'e's theorem implies that $f$ has $m_j$ zeros inside $\tau$, and these zeros must be real. Since $f'$ has no zeros in $D_j$, 
these zeros of $f$ are also simple, and $m_j \in \{ 1, 2\}$ by Rolle's theorem. 
\hfill$\Box$
\vspace{.1in}

In view of Lemma \ref{lemy2}, the hypothesis (d) may now be used for the first time, 
to separate the remainder of the proof into two cases. 
\\\\
\textbf{Case A:} \textit{assume that all but finitely many poles of $f$ have multiplicity $2$}.
\\\\
The first step in this case is the following. 

\begin{lem}
 \label{lemy3}
All but finitely many zeros of $f'$ have multiplicity $3$. 
\end{lem}
\textit{Proof.}
It is enough to take successive real zeros 
$x_{j-1} < x_j < x_{j+1} $
of $f'$ with $|x_{j-1}|$ and $|x_{j+1}|$ large, and to show that the multiplicity $n_j$ of $x_j$ is $3$. 
Since all but finitely many zeros of $f''$ are 
zeros of $f'$, Rolle's theorem implies that there exist poles $v_k$, $v_{k+1}$ of $f'$ which satisfy 
$x_{j-1} < v_k < x_j < v_{k+1} < x_{j+1}$, and these may be assumed to be the nearest poles of $f'$ to $x_j$, and to have multiplicity $3$ for $f'$.  
It then follows, using Lemmas \ref{lemboundary} and \ref{lemy1} and the argument principle, that $2 \leq n_j \leq 4$.  
On the other hand, Lemma \ref{lemy2}  and Rolle's theorem together show that  $v_k$ lies  close to, and must lie between,
a pair of real simple zeros of $f$, and the same is true of $v_{k+1}$. Thus $x_j$ lies between zeros of $f$ which are not separated by
poles of $f$, and so $x_j$ is a zero of $f'$ of  odd multiplicity, forcing $n_j = 3$. 
\hfill$\Box$
\vspace{.1in}

Now Theorem \ref{f''=0} can be applied
with $n=3$ and  $\lambda^3 =1, \lambda \neq 1$ in (\ref{f'form}), and the constants $a$ and $b$ must have zero real part. 
Hence, without loss of generality, 
$$
f'(z) = C \left( \frac{\lambda e^{ iz} - 1}{e^{ iz}-1} \right)^3 ,  
$$
and $C=1$ since $1$ is the only asymptotic value of $f'$. 
If $x$ is a pole of $f$ then, as $z \to x$,  
\begin{equation}
 \label{z1}
f'(z) \sim \frac{\mu}{(z-x)^3} , \quad f(z) \sim \frac{- \mu}{2(z-x)^2} , \quad 
\mu = \frac{(\lambda-1)^3}{i^3} = - 6 \,  {\rm Im \, }  \lambda  \in \R \setminus \{ 0 \}. 
\end{equation}

Next, let $\varepsilon_4$ be small and positive and let $U$ be the union of the discs of centre $2 \pi n$ and radius $\varepsilon_4$,
for $n \in \Z$. 
Let $m$ be an  integer with $|m|$ large, such that $m$ has the same sign as $-\mu$. Then $2 \pi m$ is a pole of $f$ and the real limit 
$\Lambda = \lim_{t \to 2 \pi m} f(t) $ exists and is infinite, with the same sign as $m$. 
Since integration shows that $f(z) \sim z$ for $z$ with $|z|$ large but $z \not \in U$, 
it follows that $\Lambda$ has the same sign as $f(2 \pi m - \varepsilon_4)$ and
$f(2 \pi m + \varepsilon_4)$. Now Rolle's theorem and the fact that $f'$ has no zeros near to $2 \pi m$ together imply that 
$f$ has no real zeros close to $2 \pi m$. But Rouch\'e's theorem gives, counting
multiplicity,  two zeros  of $f$ close to $2 \pi m$, both necessarily real, 
and this contradiction excludes Case A. 
\\\\
\textbf{Case B:} \textit{assume that all but finitely many poles of $f$ have multiplicity $1$}.
\\\\
In this case all but finitely many zeros of $f'$ have multiplicity $2$, by the argument principle. This time Theorem \ref{f''=0} may be applied
with $n=2$, and hence $\lambda= -1$, in (\ref{f'form}). This yields 
$f'(z) =  C \cot^2 (Az + B)$, with $A$, $B$, $C$  real, and the conclusion of the theorem follows easily. 
\hfill$\Box$
\vspace{.1in}

\noindent
\textit{Acknowledgement.} 
The author thanks John Rossi for invaluable discussions, and 
the referee for carefully reading  a long manuscript 
and making several very helpful suggestions and observations. 


\noindent
School of Mathematical Sciences, University of Nottingham, NG7 2RD.\\
jkl@maths.nott.ac.uk

\end{document}